\documentclass[a4paper, 12pt,leqno]{article}
\usepackage{amsmath,amsfonts,here, epsf}
\usepackage[latin1]{inputenc}
\usepackage[T1]{fontenc}
\usepackage{ae,aecompl}

\begin{document}
\newtheorem{theo}{Theorem}[section]
\newtheorem{lemme}[theo]{Lemma}
\newtheorem{cor}[theo]{Corollary}
\newtheorem{defi}[theo]{Definition}
\newtheorem{problem}[theo]{Problem}
\newtheorem{prop}[theo]{Proposition}
\newtheorem{assu}[theo]{Assumption}
\newtheorem{nontheo}[theo]{Conjectured theorem}
\newcommand{\beq}{\begin{eqnarray}}
\newcommand{\enq}{\end{eqnarray}}
\newcommand{\be}{\begin{eqnarray*}}
\newcommand{\en}{\end{eqnarray*}}
\newcommand{\Td}{\mathbb T^d}
\newcommand{\T}{\mathbb T}
\newcommand{\R}{\mathbb R}
\newcommand{\N}{\mathbb N}
\newcommand{\Rd}{\mathbb R^d}
\newcommand{\Zd}{\mathbb Z^d}
\newcommand{\Linf}{L^{\infty}}
\newcommand{\dt}{\partial_t}
\newcommand{\Dt}{\frac{d}{dt}}
\newcommand{\demi}{\frac{1}{2}}
\newcommand{\ep}{^{\epsilon}}
\newcommand{\epu}{_{\epsilon}}
\newcommand{\Dtt}{\frac{d^2}{dt^2}}
\newcommand{\vf}{\varphi}
\newcommand{\bfi}{{\mathbf \Phi}}
\newcommand{\bpsi}{{\mathbf \Psi}}
\newcommand{\bm}{{\mathbf m}}
\newcommand{\NN}{\mathbb N}
\newcommand{\RR}{\mathbb R}
\newcommand{\dx}{\partial_x}
\newcommand{\vp}{v^{\perp}}
\newcommand{\E}{{\mathbf E}}
\newcommand{\Sz}{{\mathcal{S}}}
\newcommand{\ds}{\displaystyle}
\newcommand{\fe}{f_\epsilon}
\bibliographystyle{plain}

\title{A geometric approximation to the Euler equations :
the Vlasov-Monge-Amp\`ere system}
\author{Yann Brenier, Gr\'egoire Loeper, 
\\UMR 6621, Parc Valrose, 06108 Nice, France}
\maketitle

\begin{abstract}
This paper studies the Vlasov-Monge-Amp\`ere system ($VMA$),
a fully non-linear version of the Vlasov-Poisson system ($VP$)
where the (real) Monge-Amp\`ere equation
$\det \frac{\partial^2 \Psi}{\partial x_i \partial x_j}=\rho$ 
substitutes for the usual Poisson equation.
This system can be derived as a geometric approximation of the Euler
equations of incompressible fluid mechanics in the spirit of Arnold
and Ebin. Global
existence of weak solutions and local existence of smooth solutions 
are obtained. Links between the $VMA$ system, the $VP$ system and the
Euler equations are established through rigorous asymptotic analysis.
\end{abstract}

\section{Introduction}

The classical Vlasov-Poisson ($VP$) system describes the evolution of
an electronic cloud in a neutralizing uniform background
through the following equations
\beq 
&&\frac{\partial f}{\partial t}+ \xi\cdot\nabla_x f
+\nabla_x\varphi\cdot\nabla_{\xi} f=0\label{1introprincipale}\\
&& \epsilon^2\Delta\varphi=\rho-1\label{1intropoisson},
\enq
where $f(t,x,\xi)\ge 0$
denotes the electronic density at time $t\ge 0$,
point $x\in{\mathbb  R}^d$, velocity $\xi\in{\mathbb  R}^d$
(usually $d=3$), $\rho(t,x)\ge 0$ denotes the 'macroscopic' density
\beq\label{1defrho}
\rho(t,x)=\int_{{\mathbb R}^d}f(t,x,\xi)d\xi, \label{1rho}
\enq
and $\varphi(t,x)$ denotes the electric potential at time $t$ and point
$x$ generated, through the Poisson equation (\ref{1intropoisson}),
where $\epsilon$ is a coupling constant,
by the difference between the electronic density $\rho(t,x)$
and the neutralizing background density, which is supposed to be uniform
and normalized to unity. 
Standard notations $\nabla=(\partial_1,...,\partial_d)$
and $\Delta=\partial_1^2+...+\partial_d^2$ have been used and~$\cdot$
stands for the inner product in ${\mathbb  R}^d$. 
The mathematical theory of the $VP$ system is now well understood.
In particular, existence of global
smooth solutions in three space dimensions has been proved in \cite{Pf}
(see also \cite{LP}, \cite{Sc}).
In the present paper, a fully nonlinear version of the $VP$ system is
addressed~:
\beq
&&\frac{\partial f}{\partial t}+ \xi\cdot\nabla_x f
+\nabla_x\varphi\cdot\nabla_{\xi} f=0\\
&& \det (\mathbb I + \epsilon^2 D^2 \varphi )=\rho\label{1intromonge},
\enq
where the (real) Monge-Amp\`ere equation (\ref{1intromonge})
substitutes for the Poisson equation (\ref{1intropoisson}). 
Here, $D^2\varphi(t,x)$ stands for the $d\times d$ symmetric matrix
made of all second order $x-$partial derivatives of $\varphi$,
$\mathbb I$ stands for the $d\times d$ identity matrix and $\det$
for the determinant of a square matrix.
The occurrence of the Monge-Amp\`ere equation in mathematical modeling
is not very common. Notice, however, that a very similar system 
can be found in meteorology with Hoskins' semi-geostrophic equations 
(cf. \cite{BB2}, \cite{CuG} and the included references). In a simplified
two dimensional setting, the semi-geostrophic equations read
\beq
&&\frac{\partial \rho}{\partial t}+\{\varphi,\rho\}=0
\\
&& \det (\mathbb I + \epsilon^2 D^2 \varphi )=\rho\label{semigeo},
\enq
where $\{\cdot,\cdot\}$ denotes the usual Poisson bracket.
\\
Formally, as the coupling constant $\epsilon$ is small, the $VP$ and $VMA$
equations asymptotically approach each other up to order $O(\epsilon^4)$.
Indeed, linearizing the determinant about the identity matrix leads
to
\beq
\det (\mathbb I + \epsilon^2 D^2 \varphi )=1+\epsilon^2\Delta\varphi
+O(\epsilon^4).
\enq
The formal limit, as $\epsilon=0$, reads
\beq
\label{extended euler}
&&\frac{\partial f}{\partial t}+ \xi\cdot\nabla_x f
+\nabla_x\varphi\cdot\nabla_{\xi} f=0\\
&& \rho=1\label{1introcontrainte},
\enq
where constraint (\ref{1introcontrainte}) substitutes for both the Poisson
and the Monge-Amp\`ere equations. The limit system
(\ref{extended euler},\ref{1introcontrainte}), that we call constrained
Vlasov system, can be seen as a
'kinetic' extension of the Euler equations of classical incompressible
fluid mechanics,
\beq
&&\dt v +(v\cdot\nabla) v=-\nabla p\label{1euler}\\
&&\nabla\cdot v =0\label{1div=0},
\enq
where $v(t,x)\in {\mathbb R}^d$ and $p(t,x)\in {\mathbb R}$ 
respectively are the velocity and the pressure of the fluid at time $t$
and position $x$.
Indeed, any smooth solution $(v,p)$ provides a 'monokinetic' solution to
the constrained Vlasov system (\ref{extended euler},\ref{1introcontrainte}), 
defined by
\[
f(t,x,\xi)=\delta(\xi-v(t,x)),\;\;\; \varphi=-p.
\]
Here a monokinetic solution means a delta-valued solution 
in the $\xi$ variable.
In addition, the constrained Vlasov system 
(\ref{extended euler},\ref{1introcontrainte})
turns out to be a natural
extension (or $\Gamma$ limit) 
of the Euler equations from both geometrical and variational 
reasons, as explained in section \ref{1section-euler}
\\
In a similar way, there is a monokinetic version of the $VP$ system,
the so-called (pressureless) Euler-Poisson ($EP$) system, which reads
\beq
&&\dt v +(v\cdot\nabla) v=\nabla \varphi\\
&&\dt \rho +\nabla\cdot(\rho v)=0\\
&& \epsilon^2\Delta\varphi=\rho-1.
\enq
A rigorous asymptotic analysis of the $VMA$ system
as $\epsilon \rightarrow 0$ will
be provided (sections \ref{1asymptotic} and \ref{1section-Euler-Poisson}),
in the case when the initial electronic density
\beq
f(t=0,x,\xi)=f^0(x,\xi)\label{1introprinit}
\enq
is asymptotically monokinetic, namely approaching
$\delta(\xi-v_0(x)),$
for some smooth divergence free velocity field $v_0$, as 
$\epsilon$ tends to zero.
Before this asymptotic analysis, we want to explain the
geometric origin
of the $VMA$ system.
It has been known, 
since Arnold's celebrated work (cf. \cite{AK}),
that the Euler equations (formally) describe
geodesics curves along a suitable group of volume preserving maps,
lengths being measured in the $L^2$ sense. 
We will show  (section \ref{1section-euler})
that the $VMA$ system just describes
approximate geodesics obtained through
a very natural penalty method, where $\epsilon$ stands for
the penalty parameter.
For this geometric interpretation to be valid, the Monge-Amp\`ere 
equation (\ref{1intromonge}) must be understood in the 
following weak sense:
for each fixed $t$,  $\varphi(t,\cdot)$ is the unique (up to an
additive constant) function such that 
$\Psi(x)=x^2/2 + \epsilon^2\varphi(t,x)$ is convex in $x$  and 
\beq
\forall g\in C^0(\mathbb R^d), \int_{ {\mathbb R}^d}
g(\nabla\Psi(x))\rho(t,x)dx\,=\,\int_{\Omega}g(y)dy,
\enq 
where $\Omega$ is a fixed bounded open convex set where the neutralizing
background of the electrons is assumed to be located.
(This definition is made precise in section \ref{1section-projection}.)
Notice that, by construction, $\nabla\Psi$ must be valued in the closure
of $\Omega$ and, therefore, the potential $\varphi$ enjoys the following 
property
$$
|x+\epsilon^2\nabla_x\varphi(t,x)|\le \sup_{y\in\Omega}|y|<+\infty.
$$
There is no similar bound for the electrostatic potential of the
classical $VP$ system. Thus, in some sense, the $VMA$ system can be seen
as a nonlinearly saturated version of the $VP$ system.
\\
Beyond the geometric derivation of the $VMA$ system,
our main analytic results are as follows:
\\
\begin{itemize}
\item The $VMA$ system  admits  global energy preserving weak solutions.
\item The $VMA$ system  admits  local  strong solutions in periodic domains.
\item For well prepared, nearly monokinetic initial data, the solutions of the $VMA$ system 
converge when $\epsilon$ goes to 0 to those of the Euler equations.
\item In this asymptotic,  the $EP$ system is a higher
order approximation of the $VMA$  system.
\end{itemize}
\bigskip
The paper is organized as follows: in section 2,
we first recall the geometric nature of the Euler equations, then we explain
why the constrained Vlasov system 
(\ref{extended euler},\ref{1introcontrainte}) is a natural
extension of the Euler equations from a variational point of view,
finally we introduce the concept of approximate geodesics 
for volume preserving maps, 
and derive the $VMA$ system.
Section 3 is devoted to the proof of existence of global 
energy preserving weak solutions. In section 4,
we prove existence of local strong solutions, in the case of a periodic
domain. Finally, in section 5, we study the asymptotic behavior of the
$VMA$ system as $\epsilon$ goes to 0.

\section{The geometric origin of the Vlasov-Monge-Amp\`ere system}
\label{1section-euler}
\subsection{The Euler equations}
The motion of an incompressible fluid in a domain $\Omega\subset\Rd$ 
is classically described by the
 Euler equations $(E)$:
\beq
&&\dt v +(v\cdot\nabla) v=-\nabla p\\
&&\nabla\cdot v =0,
\enq
with $t\in \R$, $x\in \Omega$, where $v=v(t,x)$ 
stands for the velocity field and $p=p(t,x)$ for the scalar pressure field.
These equations have a nice geometrical interpretation
going back to Arnold (see \cite{AK}).
Introducing $G(\Omega)$ the group of  
all volume preserving diffeomorphisms of $\Omega$
with jacobian determinant equal to 1, and measuring lengths in the $L^2$ sense,
we may define
(at least formally) geodesic curves along $G(\Omega)$. It turns out
that the Euler equations just describe these curves. 
For the same reasons, the Euler equations can be seen as the optimality equations
for the corresponding minimization problem: given two maps chosen in
$G(\Omega)$, find an $L^2-$shortest path between them along 
$G(\Omega)$. It was shown by Shnirelman \cite{Shn1} (see also
\cite{AK} and \cite{Shn2}) that, in the case
when $\Omega$ is the unit cube in $\R^3$, there are many maps for which there
are no such shortest paths. 
Beyond this negative result, \cite{Br5} established that 
minimizing paths are more appropriately described by 
doubly stochastic measures. These measures (also called polymorphisms)
generalize volume preserving maps in the following way:
a doubly stochastic measure $\mu(dx,dy)$ is a (Borel) probability measure on 
$\Omega\times\Omega$ with two projections on each copy of $\Omega$ both equal
to the (normalized) Lebesgue measure. It is known  -see \cite{Neretin},
for instance-
that any such $\mu$ can be weakly approximated by a sequence 
$\mu_n(dx,dy)=\delta(x-g_n(y))dy$ where each $g_n$ is a volume preserving
map of $\Omega$.
In \cite{Br5} it was shown that, in the case considered by Shnirelman for which
there is no classical shortest path, minimizing paths along $G(\Omega)$
converge to paths of doubly stochastic measures $t\rightarrow \mu(t;dx,dy)$
governed by the following extension of the Euler equations
\beq
&&\dt \mu+\nabla_x\cdot(\mu v)=0,\label{varia1}\\
&&\dt (v\mu) +\nabla_x\cdot(\mu v\otimes v)+\mu\nabla_x p=0,\label{varia2}
\enq
where $v=v(t;x,y)$ and $p=p(t,x)$ can be respectively seen as the
velocity field and the pressure field attached to $\mu$.
(Notice that the velocity field $v$ generally depends on the extra variable $y$
and is not a classical but rather a multivalued
velocity field.)
These equations are just a reformulation of
the constrained Vlasov system (\ref{extended euler},\ref{1introcontrainte}).
Indeed, it can be checked, under appropriate regularity assumptions, 
that the kinetic measure $f$ defined by
\beq\label{defmuf}
f(t;dx,d\xi)=\int_{y\in\Omega} \delta(\xi-v(t;x,y))\mu(t;dx,dy)
\enq
solves (\ref{extended euler},\ref{1introcontrainte})
when $(\mu,v,p)$ solves (\ref{varia1},\ref{varia2}).
Thus we conclude 
that the constrained Vlasov system 
(\ref{extended euler},\ref{1introcontrainte})
is a natural variational extension of the Euler equations.

\subsection{Approximate geodesics }\label{1section-apxgeo}
A general strategy to define approximate geodesics along a manifold $M$
(in our case $M=G(\Omega)$) embedded in a Hilbert space $H$ 
(here $H=L^2(\Omega,\Rd)$)
is to introduce a penalty parameter $\epsilon>0$
and the following $unconstrained$ dynamical system in $H$
\beq
\partial_{tt}{X} +\frac{1}{2\epsilon^2}\nabla_X\left(d^2(X,M))\right)
=0.
\label{1apxgeo}
\enq
In this equation,
the unknown $t\rightarrow X(t)$ is a curve in $H$, 
$d(X,M)$ is the distance (in $H$) of $X$ to the manifold
$M$, i.e. in our case as $M=G(\Omega)$,
\beq
d(X,G(\Omega))=\inf_{g\in G(\Omega)}\|X-g\|_H,\label{1projection}
\enq
and, finally, $\nabla_X$ denotes the
gradient operator in $H$.
This penalty approach 
has been used for the Euler equations 
by the first author in \cite{Br3}. 
It is similar-but not identical-
to Ebin's slightly compressible flow theory \cite{Eb}, and is
a natural extension of the theory of 
constrained finite dimensional mechanical systems \cite{RU}. 
The penalized system is formally hamiltonian in variables
$(X,\partial_t X)$ with 
Hamiltonian (or energy) given by:
\be
E=\frac{1}{2}\|\partial_t X\|_H^2 + \frac{1}{2\epsilon^2}d^2(X,G(\Omega)).
\en
(Multiplying equation (\ref{1apxgeo}) by $\partial_t X$, we formally get
that the energy is conserved.)
Therefore it is plausible that the map $X(t)$ will remain close to $G(\Omega)$
if properly initialized at $t=0$.
A formal computation shows that, given a point $X$ for which there is
a unique closest point $\pi_X$ to $X$ in the $H$ closure
of $G(\Omega)$, 
we have:
\beq
\nabla_X\left(d(X,G)\right)=\frac{1}{d(X,G)}(X-\pi_X).\label{1grad}
\enq
Thus the equation (\ref{1apxgeo}) formally becomes: 
\beq
\label{1apxgeo2}
\partial_{tt} X + \frac{1}{\epsilon^2}(X-\pi_X)=0.
\enq
To understand why solutions to such a system may approach
geodesics along $G(\Omega)$ as 
$\epsilon$ goes to 0, just recall that, in the simple framework of a 
surface $S$ embedded in the 3 dimensional Euclidean space,
a  geodesic $t\rightarrow s(t)$  along $S$ 
is characterized by the fact that
for every $t$, the plane defined by $\{\dot s(t), \ddot s(t)\}$ is 
orthogonal to $S$.
In our case, $\partial_{tt} X(t)$ is nearly 
orthogonal to $G(\Omega)$ thanks to 
(\ref{1apxgeo2}), meanwhile $X(t)$ remains close to $G(\Omega)$.
\\
The approximate geodesic equation was introduced in \cite{Br3}
in order to allow a spatial approximation of $G(\Omega)$
by the group of permutations of $N$ points $A_j$ chosen to form a discrete
grid on $\Omega$. On such a discrete group, the concept of geodesics
becomes unclear meanwhile approximate geodesics still make sense.
They can be interpreted as trajectories of a cloud of $N$ particles $X_i$
moving in the Euclidean space ${\mathbb R}^{dN}$, which substitutes
for $H$. These particles solve the following coupled system of
harmonic oscillators
\[
\epsilon^2\frac{d^2 X_i}{dt^2}+ X_i-A_{\sigma_i}=0,
\]
where $\sigma$ is a time dependent permutation minimizing,
at each fixed time $t$,
$\Sigma \left|X_i -A_{\sigma(i)}\right|^2$ among all other permutations 
of the first $N$ integers.
The  convergence of this discrete model to the incompressible 
Euler equations for well prepared initial data was proved in \cite{Br3}.
In order  to  study the continuous version (\ref{1apxgeo2}),
a specific study of the projection problem (\ref{1projection}) is needed.

\subsection{The polar decomposition Theorem}\label{1section-projection}

Let us first recall a general measure theoretic definition:
\begin{defi}
Let $A$ and $B$ be two topological spaces, let  
$\rho$ be a Borel finite 
measure of $A$ and  $X$ a Borel map  $A\rightarrow B$, 
we call the push-forward of $\rho$ by $X$ and note $X\#d\rho$ the Borel
measure $\eta$ on $B$ defined by
\be
\forall f\in C^0(B),\;\int_{B}f(y) d\eta(y) = \int_{A}f(X(x))d\rho(x).
\en 
\end{defi}
Let us now consider the case of a bounded open subset $\Omega$ of
the Euclidean space $\Rd$ equipped with the Lebesgue 
measure that we denote $dx$.
We say that a Borel map $s: \overline{\Omega}\rightarrow \overline{\Omega}$ is volume
(or Lebesgue measure) preserving 
if $s\#dx=dx$, i.e. if for all $g\in C^0(\overline{\Omega})$ one has 
$\int_{\Omega}  g(x)dx =\int_{\Omega}  g(s(x))dx,$
or equivalently, for any Borel subset $B$ of $\overline{\Omega}$ one has $|s^{-1}(B)|=|B|$.
The set of all measure preserving maps of $\Omega$ is 
a closed subset of the Hilbert space $H=L^2(\Omega,\Rd)$
and will be denoted by $S(\Omega)$.
Notice that $S(\Omega)$ is only a semi-group for
the composition rule and contains the group
of volume preserving diffeomorphisms $G(\Omega)$.
It is known \cite{Ne} that, at least in the case when $\Omega$ is convex
and $d\geq 2$, $S(\Omega)$ is exactly the closure of $G(\Omega)$ 
in $L^2(\Omega,\Rd)$, which implies $d(.,G(\Omega))=d(.,S(\Omega)).$
\\
The polar decomposition Theorem for maps \cite{Br1} (extended to
Riemannian manifolds in \cite{Mc2})
will be crucial for
our analysis of the $VMA$ system:

\begin{theo}\label{1brenier}
Let $\Omega$ be a bounded convex open subset of $\Rd$, 
let $X\in L^2(\Omega;\Rd)$ 
and $\rho_X=X\#dx,$ where $dx$ is the Lebesgue measure on $\Omega.$ 
Assume $\rho_X$ to be a Lebesgue integrable function,
or, equivalently, $X$ to satisfy the non-degeneracy condition:
\beq
\forall E\subset \Rd \mbox{ Borel },\, |E|=0\Rightarrow |X^{-1}(E)|=0.
\label{1nondeg}
\enq
Then there exists a unique pair $(\nabla\Phi_X, \pi_X)$ 
where $\Phi_X$ is a convex function and 
$\pi_X\in S(\Omega)$,  such that
\beq
X=\nabla\Phi_X\circ \pi_X.\label{1polaire}
\enq
In this 'polar decomposition',
$\pi_X$ is also characterized as the unique
closest point to $X$ on $S(\Omega)$ in the $L^2$ sense
and
$\Phi_X$ is characterized to be (up to an additive constant) 
the unique convex function 
on $\Omega$ satisfying
\beq
\int_{\Rd}g(x)d\rho_X=\int_{\Omega}g(X(y))dy=\int_{\Omega}g(\nabla\Phi_X(y))dy\label{1trans1},
\enq 
for any $ g \in C^0(\Rd)$ such that $\left|g(x)\right|\leq C(1+|x|^2)$.
\\
In addition, the Legendre-Fenchel transform $\Psi_X$ of $\Phi_X$ defined by
\beq \label{1defLegendre}
\Psi_X(x)=\sup_{y\in \Omega}\{x\cdot y -\Phi_X(y)\}
\enq
is Lipschitz continuous on $\Rd$, with Lipschitz constant bounded by
$\sup_{x\in\Omega}|x|$ and has the following properties~:
\\
$\nabla\Psi_X(x)\in\Omega$ holds true for $\rho_X$ a.e. $x$, 
\beq
\int_{\Rd}g(\nabla\Psi_X)\rho_X(x)dx=\int_{\Omega}g(\nabla\Psi_X(X(x)))dx=\int_{\Omega}g(x)dx\label{1trans2}
\enq 
for any  $g \in C^0(\overline{\Omega})$, and
\beq
&&\nabla\Phi_X(\nabla\Psi_X(x))=x \;\;\;\rho_X(x)dx \;a.e,\label{1inverse1}\\
&&\nabla\Psi_X(\nabla\Phi_X(y))=y\;\;\;dy \;a.e,\label{1inverse2}\\
&&\pi_X(y)=\nabla\Psi_X(X (y))\;\;\;dy \;a.e.\label{1inverse3}
\enq

\end{theo}
We make here several remarks on Theorem \ref{1brenier}:

\bigskip
\noindent 
{\bf Link with the Monge-Amp\`ere equation} 
We can interpret (\ref{1trans1}) 
as a weak version of the Monge-Amp\`ere equation:
\be
\rho_X(\nabla \Phi)\det D^2\Phi=1
\en and (\ref{1trans2}) can be seen as a weak version of another Monge-Amp\`ere equation:
\be
&&\det D^2\Psi=\rho_X\\
&& \nabla\Psi \mbox{ maps supp}(\rho_X)\mbox{ in }\Omega.
\en

\bigskip
\noindent 
The pair $(\Phi_X, \Psi_X)$ depends in fact only of $\Omega$ and the measure $\rho_X=X\#dx,$ 
and if condition (\ref{1nondeg}) fails, 
then existence and uniqueness of the projection $\pi_X$ 
may fail, but existence and uniqueness of $\nabla\Phi_X$ remain true. 

Theorem \ref{1brenier} and the subsequent remarks allow us to introduce the following
 notation that will be used throughout the paper:

\begin{defi}\label{1MaOmegarho}
Let $\Omega$ be a fixed bounded convex open set of $\Rd$, let $\rho$ be a positive measure
 on $\Rd$ of total mass $|\Omega|$, absolutely continuous w.r.t the Lebesgue measure and such that $\int (1+|x|^2) d\rho(x) < +\infty$. We call 
 $\Phi[\Omega,\rho]$, or, in short, $\Phi[\rho]$,   
the unique up to a constant  convex function on $\Omega$  satisfying 
\beq
\forall g\in C^0(\mathbb R^d)\cap L^1(d\rho{}), \int_{ {\mathbb R}^d}g(x)d\rho(x)\,=\,\int_{\Omega}g(\nabla\Phi[\Omega,\rho](y))dy.
\enq We call $\Psi[\Omega,\rho]$  its Legendre-Fenchel transform
 satisfying 
\beq
\forall g\in C^0(\mathbb R^d)\cap L^1(\Omega, dy), \int_{ {\mathbb R}^d}g(\nabla\Psi[\Omega,\rho](x))d\rho(x)\,=\,\int_{\Omega}g(y)dy.
\enq 
\end{defi}
If no confusion is possible we may write $\Phi$ (resp. $\Psi$)
 instead of $\Phi[\Omega,\rho]$ (resp. $\Psi[\Omega,\rho]).$

\bigskip
\noindent
We will use some additional results from
\cite{Br1}. The first one establishes 
the continuity of the polar decomposition:
\begin{theo}\label{1convergence}
Let $\rho$ be a Lebesgue integrable positive measure on $\Rd$, with total mass $\Omega$, such that $\int (1+|x|^2) d\rho < +\infty$.
Let $\rho_n$ be a sequence of Lebesgue integrable positive measures on $\Rd$, with total mass $\Omega$, such that 
$\forall n$, \mbox{$\int (1+|x|^2) d\rho_n < +\infty$}.
Let $\Phi_n=\Phi[\Omega,\rho_n]$ and 
$\Psi_n=\Psi[\Omega,\rho_n]$ be as in Definition \ref{1MaOmegarho}. 
If for any $f\in C^0(\Rd)$ such that $|f(x)|\leq C(1+|x|^2),$  $\int f \ d\rho_n$ converges to  $\int f \ d\rho $,
  then 
\begin{itemize}
\item $\Phi_n$ converges to $\Phi[\Omega,\rho]$
 uniformly on each compact set of $\Omega$ and strongly in $W^{1,1}(\Omega),$
\item $\Psi_n$ converges to $\Psi[\Omega,\rho]$ uniformly on each compact set of $\Rd$ and strongly in $W^{1,1}(K)$ for every $K$ compact in $\Rd$.
\end{itemize} 
\end{theo} 
The second one provides a 'dual' definition of the distance between
a map $X$ and the semi-group $S(\Omega)$:
\begin{theo}\label{1distance}
Let $X\in L^2(\Omega;\Rd)$ 
and $\rho=X\#dx,$ where $dx$ is the Lebesgue measure on $\Omega.$ 
Assume $\rho$ to be a Lebesgue integrable function.
Then
$$
\frac{1}{2}d^2(X,S(\Omega))
=\int \left(|x|^2/2-\Psi[\Omega,\rho](x)\right)\rho(x)dx 
+\int_{\Omega} \left(|y|^2/2-\Phi[\Omega,\rho](y)\right)dy\\
$$
$$
=\sup_{u,v}
\int \left(|x|^2/2-u(x)\right)\rho(x)dx +\int_{\Omega} \left(|y|^2/2-
v(y)\right)dy,
$$
where the supremum if performed over all pairs $(u,v)$ of continuous
functions on $\Rd$ such that $u(x)+v(y)\ge x\cdot y$ pointwise.
\end{theo} 

\subsection{The Vlasov-Monge-Amp\`ere system}

Let us now derive the $VMA$ system as the  kinetic formulation
of the approximate
geodesic equation (\ref{1apxgeo2}). 
First, from the polar decomposition Theorem \ref{1brenier}, equation 
(\ref{1apxgeo2}) reads
\beq
\partial_{tt} X(t,x) = \nabla\varphi(t,X(t,x)),\label{1odeapxgeo}
\enq
where 
\beq
\label{1varphi}
\nabla\varphi(t,x)=\frac{\nabla\Psi[\Omega,\rho(t,\cdot)](x)-x}{\epsilon^2}
\enq
and $\Psi[\Omega,\rho]$ is as in Definition (\ref{1MaOmegarho}).
This means that $\nabla\varphi$ satisfies (\ref{1intromonge}) in
a weak form with the additional
condition that the range of $x\rightarrow x+\epsilon^2\nabla\varphi(t,x)$ is
contained in $\overline{\Omega}$.

Next, let $f^0\ge 0$ be a given initial density function, that
we assume to be in $L^{\infty}(\Rd\times\Rd)$,  compactly supported and 
satisfying the compatibility condition
\beq
\int  f^0(x,\xi)dxd\xi=|\Omega|.\label{1jauge}
\enq
For each $t\ge 0$, let us define $(x,\xi)\rightarrow f(t,x,\xi)$ to be 
$f^0$ pushed forward 
by the following ODE 
\beq
&&\partial_t X(t,x,\xi)=\Xi(t,x,\xi)\\
&&\partial_{t} \Xi(t,x,\xi) = (\nabla\varphi)(X(t,x,\xi))\label{1od1}\\
&& (X,\Xi)(t=0,x,\xi)=(x,\xi)\label{1od2}.
\enq
Then $f$ satisfies the following kinetic (or Liouville) equation
\beq
&&\frac{\partial f}{\partial t}+\nabla_x\cdot\left(\xi f\right)
+\nabla_{\xi}\cdot\left(\nabla\varphi f\right)=0\label{1principale}\\
&&f(0,\cdot,\cdot)=f^0,\label{1prinit}
\enq
which must be understood in the following weak sense
\beq
&&\forall g \in C^{\infty}_c([0,+\infty)\times \Rd \times \Rd),\nonumber\\
&&\int_0^\infty dt
\int_{\Rd\times\Rd} \left(\frac{\partial g} {\partial t} + \xi\cdot \nabla_x g 
+  \nabla\varphi\cdot\nabla_{\xi} g \right )  f  dx d\xi\nonumber\\
&& = -\int_{\Rd\times\Rd}
f_0(x,\xi)g(t=0,x,\xi)dx d\xi.
\enq
This linear Liouville equation is
nonlinearly coupled to equation (\ref{1varphi}),
where $\rho$ is linked to $f$ by equation (\ref{1rho}). 
Finally, we have defined, through 
(\ref{1varphi},\ref{1principale},\ref{1prinit}), the weak formulation
of the $VMA$ initial value problem.

\noindent
The energy of the system is defined by
\beq
E(t)=&&\frac{1}{2}\int_{\Rd\times\Rd}f(t,x,\xi)|\xi|^2dxd\xi\nonumber\\
+&&\frac{1}{2\epsilon^2}\int_{\Rd}\rho(t,x)\left|\nabla\Psi[\Omega,\rho](t,x)-x\right|^2dx\label{1engdef}.
\enq


\section{Existence of global renormalized weak solutions}
The main result of this section is as follows:
\begin{theo}\label{1globalweak}
Let $(x,\xi)\rightarrow f^0(x,\xi)\ge 0$ be in 
$L^{\infty}(\Rd\times\Rd)$, with compact support in both $x$  
and $\xi$, satisfying condition (\ref{1jauge}).
\\
Then the $VMA$ system (\ref{1varphi},\ref{1principale},\ref{1prinit}) 
admits a global
weak solution $(f,\rho,\Psi)$, with $f\in L^{\infty}(\R^+\times \Rd\times\Rd)$
and $(\rho, \nabla\psi) \in \Linf([0,T]\times\Rd)$ for all $T>0$. In addition,
each such weak solution enjoys the following properties:
\begin{itemize}
\item $f$ is a continuous function of $t$, valued in
$L^p(\Rd\times\Rd)$, for every $1\leq p <\infty$
\item the density $\rho$ is a continuous function of $t$, valued in
$L^p(\Rd)$, for every $1\leq p <\infty$,
\item the support of $f(t,\cdot,\cdot)$ in $(x,\xi)$ is
compact, with a diameter growing no more than linearly in $t$.
\item the total energy defined by (\ref{1engdef}) is conserved,
\item the 'renormalization' property (in the sense of \cite{DL})
$$\frac{\partial g(f)}{\partial t}+\nabla_x\cdot\left(\xi g(f)\right)
+\nabla_{\xi}\cdot\left( \nabla\varphi g(f)\right)=0$$ 
holds true for  all $g\in C^1(\mathbb R)$,
\item the trajectories of (\ref{1od1},\ref{1od2}) are uniquely defined
for almost every initial condition $(x,\xi)$,
\item $t\rightarrow f(t,\cdot,\cdot)$ is just $f^0$ pushed forward 
along the trajectories of (\ref{1od1},\ref{1od2}).
\end{itemize}
\end{theo}
\noindent
{\bf Proof of Theorem \ref{1globalweak}:}

\noindent
We build a sequence of approximate 
solutions $(f_h,\Psi_h)_{h>0}$ by time discretization and let the
time step $h$ go to zero.
To handle the limiting process,
the non-linear terms will be treated with the help
of  Theorem \ref{1convergence}. 
More precisely if one can extract a subsequence
 such that, for every $t$,
 $f_h(t,\cdot,\cdot)$ converges weakly, then we can deduce from Theorem 
\ref{1convergence} that the corresponding sequence $\nabla\Psi_h(t,\cdot)$ 
will converge strongly, and this will allow us pass to
 the limit in the nonlinear term. 
 \subsection{Construction of a sequence of approximate solutions}
We consider $\eta\in C^{\infty}_c(\Rd)$ such that $\eta\geq 0$, $\int_{\Rd}\eta =1$ and
$\eta_h=\frac{1}{h^d}\eta(\frac{\cdot}{h})$. We then seek approximate solutions as 
solutions of the approximate problem
\beq
&& \frac{\partial f_h}{\partial t}+\xi\cdot\nabla_xf_h
+\frac{\nabla\Psi_h(x)-x}{\epsilon^2}\cdot\nabla_{\xi} f_h=0\\
&&f_h(0,x,\xi)=f_h^0(x,\xi)=f_0*_{x,\xi} \eta_h\otimes\eta_h\\
&&\Psi_h(t)=\eta_h * \Psi[\Omega,\rho(t=nh)] \mbox{ for } t\in [nh,(n+1)h[.
\enq
$\nabla\Psi_h$ being a smooth function of space this regularized equation admits a unique solution 
that one builds by the method of characteristics.
Since the flow is
 divergence-free in the phase space, the solution $f_h$ satisfies
\beq \label{1consnorme}
\forall p \in [1,+\infty],\;\|f_h(t)\|_{L^p(\Rd\times\Rd)}=\|f_h(0)\|_{L^p(\Rd\times\Rd)}.
\enq
By construction (through Theorem \ref{1brenier}),
$\nabla\Psi_h$ is valued in the convex bounded set $\overline{\Omega}$.
Suppose that $f^0(x,\xi)$ vanishes outside of the set
 $\{x^2+\epsilon^2\xi^2\leq C^2\}$ for some constant $C>0$ fixed and 
denote $R=\sup_{y\in \Omega}|y|$. Then we have
\begin{lemme}\label{1suppcpct}
$\forall t\geq 0,\; f_h(t,\cdot,\cdot)$ is supported in   
$\{\sqrt{x^2+\epsilon^2\xi^2}\leq C+ Rt/\epsilon\}$.
\end{lemme}
{\bf Proof :}  
We just write
$$\epsilon^2 \partial_{tt}{X} + X = \nabla\Psi_h(X)$$
in complex notation
$-i\epsilon\partial_t{Z}+Z=F$, where $Z=X+i\epsilon\partial_t{X}$ and $F=\nabla\Psi_h(X)$,
which is bounded by $R$. This leads to
$$Z(t)=Z(0)\exp(-it/\epsilon)+i\epsilon^{-1}
\int_0^t \exp(-i(t-s)/\epsilon)F(s)ds$$
ant the desired bound easily follows.
Notice here a sharp contrast with the classical $VP$ system, for which
the $\xi-$support of the solutions 
cannot be controlled so easily (except in the one dimensional case).
\hfill $\Box$

\bigskip
\noindent
{\bf Convergence of the sequence of approximate solutions}
\\
\noindent
Using (\ref{1consnorme}) and Lemma \ref{1suppcpct}  there exists, for any $1<p<\infty$, 
up to the extraction of a subsequence, $f \in L^p([0,T]\times \Rd \times \Rd)$
such that $f_h $ converges weakly to $f$ as $h\rightarrow 0$.

It remains to show that the product $f_h \nabla\Psi_h$ converges to the good limit.
For this we need strong convergence of $\nabla\Psi_h.$ We already know that 
$\nabla\Psi_h  \in L^{\infty}([0,T]\times \Rd).$
We claim that for all $t>0$, $\nabla\Psi_h(t,\cdot)$ 
converges strongly to $\nabla\Psi(t,\cdot)$
in $L^q_{loc}(\Rd)$, $\forall q\in [1,+\infty[$. Indeed, such a strong convergence 
of $\nabla\Psi_h$ follows from Theorem \ref{1convergence} provided 
that we have for all $t>0$,
\beq
\int_{{\mathbb R}^d}g(x)\rho_h(t,x)dx\rightarrow
\int_{{\mathbb R}^d}g(x)\rho(t,x)dx,
\label{1MA}
\enq
for any  $g\in C^0({\mathbb R}^d)$ such that $\int (1+|x|^2)g(x)dx<+\infty$.
Note first that from Lemma \ref{1suppcpct}, we can restrict ourselves here to test functions $g$ that are compactly supported.
Then we show that  the sequence $\rho_h$ is  relatively compact in 
 $C([0,T],L^p(\Rd)-w).$ This is done by the following lemma:
\begin{lemme}  For all $T>0$, for all $p$ with $1\leq p< \infty$ the sequence
 $f_h$ (resp. $\rho_h$) satisfies
\begin{itemize}
\item $f_h$ (resp. $\rho_h$) is a bounded sequence in $\Linf([0,T];L^p(\Rd\times\Rd))$ 
(resp. in $\Linf([0,T];L^p(\Rd))$,
\item $\dt f_h$ (resp. $\dt \rho_h$) is a bounded sequence in $\Linf([0,T];W^{-1,p}(\Rd\times\Rd)))$,
(resp. in $\Linf([0,T];W^{-1,p}(\Rd))$,
\end{itemize}
and  one can extract from $f_h$  (resp. from $\rho_h$) a subsequence converging in $C([0,T],L^p(\Rd\times\Rd)-w)$ (resp. in $C([0,T],L^p(\Rd)-w)$).
\end{lemme} 
{\bf Proof:} the first point uses equation (\ref{1consnorme}) and Lemma \ref{1suppcpct}.
 The second point uses equation (\ref{1principale}) and the identity:
$$\dt \rho_h=-\nabla_x \cdot \int_{\Rd} \xi f_h d\xi,$$
with the fact that the $f_h$ are uniformly compactly supported in $x$ and $\xi$ (Lemma \ref{1suppcpct}); 
the last point is a classical
 result of functional analysis (see \cite{Li} for example). $\hfill$ $\Box$

\bigskip
\noindent
This lemma and Lemma \ref{1suppcpct} yield (\ref{1MA}). Then using Theorem \ref{1convergence},
 with $\rho$ the limit of a subsequence of $\rho_h$, we have convergence
 of the sequence $\nabla\Psi_h$ to  $\nabla\Psi[\Omega,\rho]$ in $C([0,T],L^p(\Rd))$.
We have extracted a subsequence $f_h$ such that
\begin{itemize}
\item $f_h$ converges in $C([0,T],L^p(\Rd\times\Rd)-w)$ for every $1\leq p <\infty$.
\item $\rho_h$ converges in $C([0,T],L^p(\Rd)-w)$ for every $1\leq p <\infty$.
\item $\nabla\Psi_h(t,\cdot)$ converges in $L^p(\Rd)$ for every $t$ and
for every $1\leq p <\infty$\label{1Psiconv}.
\end{itemize}
Thus the limit $(f,\nabla\Psi)$ satisfies equations (\ref{1principale}-\ref{1prinit})
 and the first part of Theorem \ref{1globalweak} is proved. 


\subsection{Conservation of energy}\label{1energysection}
We now give a rigorous proof of the conservation of energy 
following an argument going back to F. Otto (in an unpublished
work on the semi-geostrophic equations).
We recall the definition of the energy as
\be
E(t)=\frac{1}{2}\int_{\Rd\times\Rd}f(t,x,\xi)|\xi|^2dxd\xi+\frac{1}{2\epsilon^2}\int_{\mathbb R^d}\rho(t,x)|\nabla\Psi(t,x)-x|^2dx.
\en We call the first term the kinetic energy $E_c$  and the  second term,
multiplied by $\epsilon^2$, the (normalized) potential energy $E_p.$
We have
\begin{prop}
Let $f$ be any solution of (\ref{1principale}) such that on every interval $[0,T]$, $f(t,\cdot,\cdot)$ is uniformly compactly supported in $|x|, |\xi|\leq R(T)$ for some function $R(T)$. Then the energy of the solution $f$ is conserved. 
\end{prop}
{\bf Proof:}
From Theorem \ref{1distance},
we know that
$$
E_p(t)=\int \left(|x|^2/2-\Psi(t,x)\right)\rho(t,x)dx 
+\int_{\Omega} \left(|y|^2/2-\Phi(t,y)\right)dy
$$
$$
=\sup_{u,v}
\int \left(|x|^2/2-u(x)\right)\rho(t,x)dx +\int_{\Omega} \left(|y|^2/2-
v(y)\right)dy,
$$
where the supremum if performed over all pairs $(u,v)$ of continuous
functions on $\Rd$ such that $u(x)+v(y)\ge x\cdot y$ pointwise.
Thus for each $t, t_0 \in \mathbb R+$, we have
\be
E_p(t)&&\geq \int \left(|x|^2/2-\Psi(t_0,x)\right)\rho(t,x)dx +\int_{\Omega} \left(|y|^2/2-\Phi(t_0,y)\right)dy,
\en
and this implies
\begin{eqnarray*}
E_p(t)-E_p(t_0)&&\geq \int_{\mathbb R^d} \left(|x|^2/2-\Psi(t_0,x)\right)(\rho(t,x)-\rho(t_0,x))dx\\
&&= \int_{t_0}^{t}\int_{\mathbb R^d}\partial_t\rho(s,x)\left(|x|^2/2-\Psi(t_0,x)\right)dxds\\
&&= \int_{t_0}^{t}\int_{\R^d\times \Rd}\xi f(s,x,\xi)\left(x-\nabla\Psi(t_0,x)\right)dxd\xi ds.
\end{eqnarray*}
Notice that the product in the second line is licit
 since $\partial_t\rho$ is
 in $W^{-1,p}$ for any $1\leq p < \infty$, 
$f(t,\cdot,\cdot)$ and therefore $\rho(t,\cdot)$ are compactly supported in space uniformly on $[0,T]$,
and $\Psi-|x|^2/2$ is in $W^{1,\infty}_{loc}$.
Exchanging $t_0$ and $t$ we would have found
\be
E_p(t_0)-E_p(t)\geq\int_{t}^{t_0}\int_{\R^d\times\Rd}\xi f(s,x,\xi)\left(x-\nabla\Psi(t,x)\right)dxd\xi ds,
\en 
moreover we have for the kinetic energy
\be
\epsilon^2(E_c(t)-E_c(t_0))= 
\int_{t_0}^{t}\int_{\R^d\times\Rd} \xi f(t,x,\xi)\cdot(\nabla\Psi(s,x)-x)dxd\xi ds.
\en
Dividing by $t-t_0,t> t_0$ we find
\be
&&\epsilon^2\frac{E(t)-E(t_0)}{t-t_0}\\
&&\geq\frac{1}{t-t_0}\int_{t_0}^{t} \int_{\R^d\times\Rd}\xi f(s,x,\xi)\cdot(\nabla\Psi(s,x)-\nabla\Psi(t_0,x))dxd\xi ds
\en
and
\be
&&\epsilon^2\frac{E(t)-E(t_0)}{t-t_0}\\
&&\leq\frac{1}{t-t_0}\int_{t_0}^{t}\int_{\Rd\times\Rd} \xi f(s,x,\xi)\cdot(\nabla\Psi(t,x)-\nabla\Psi(s,x))dxd\xi ds.
\en
We know from \ref{1Psiconv} that $\nabla\Psi(t,.)$ converges strongly in $L^p_{loc}(\Rd), 1\leq p<\infty$  to $\nabla\Psi(t_0,.)$ as $t$ goes to $t_0$, and so the  right hand sides of the above inequalities converges to 0 
and we conclude that

$$\lim_{t>t_0}\frac{E(t)-E(t_0)}{t-t_0}= 0.$$
We could take $t<t_0$ and find the same result. Finally we conclude that 
\be
\frac{dE}{dt}\equiv 0.
\en 
$\hfill$ $\Box$


\subsection{Renormalized solutions and existence of characteristics}
The study of renormalized solutions for transport equations 
has been introduced in \cite{DL} for 
vector fields in  $W^{1,1}$ with  bounded divergence.  
These results have been extended by Bouchut \cite{Bo}
to the case of Vlasov-type 
 equations with acceleration field in $BV$ 
(A recent result of L. Ambrosio, \cite{Am}, has extended the existence of renormalized solutions to transport equations with vector fields in $BV$ and with bounded divergence).
The fact that solutions of (\ref{1principale}, \ref{1prinit}) are renormalized solutions is an immediate consequence of the following theorem:

\begin{theo}\label{1Bouchut} $(${\rm F. Bouchut}$)$
\\
Let $f \in L^{\infty}(]0,T[,L^{\infty}_{loc}( \Rd \times \Rd))$ satisfy
\begin{eqnarray*}
\frac{\partial f}{\partial t}+\nabla_x\cdot\left(\xi f\right)
+\nabla_{\xi}\cdot\left( E(t,x)f\right)=0,
\end{eqnarray*}
with  $E(t,x)\in L^1(]0,T[;L^1_{loc}(\Rd))\cap L^{1}(]0,T[;BV_{loc}(\Rd))$,
\\
then, for 
any $g\in C^1(\R)$, 
$$
\frac{\partial g(f)}{\partial t}+\nabla_x\cdot\left(\xi g(f)\right)
+\nabla_{\xi}\cdot\left( E(t,x)g(f)\right)=0,
$$
and  for every $1\leq p <\infty$, $f$ belongs to $C(]0,T[, L^p_{loc}(\Rd\times\Rd))$.

\end{theo}

In our case the BV bound on the
 acceleration $\nabla\Psi$ 
is a direct consequence of the fact that $\Psi$ is a globally 
Lipschitz convex function.
This result implies the strong time continuity results for $f$ and  $\rho$ in Theorem \ref{1globalweak}.
Finally, as in \cite{DL}, it can be deduced from
the renormalization property that
\\
1) for almost every initial condition $(x,\xi)$,
there is a unique trajectory solving (\ref{1od1},\ref{1od2}),
\\
2) $t\rightarrow f(t)$ is just $f^0$ pushed forward along these trajectories.

A complete proof is given in appendix.

Remark: From the renormalization property it follows that,
once the potential $\Psi(t,x)$ is known, there exists a unique solution to 
(\ref{1principale})
 in $\Linf_{t,x,\xi}$. Of course, this does not imply at all 
the uniqueness of weak solutions to the
 Vlasov-Monge-Amp\`ere system!
This paragraph ends  the proof of Theorem \ref{1globalweak}.


\section{Strong solutions}
In this section we show existence of strong solutions over a finite time intervall. To do so, we need regularity estimates for solutions of 
Monge-Amp\`ere equation. We will get rid of the difficulties that may arise at the free boundary of the set $\{\rho>0\}$ by considering the periodic case. 
Note that for the Vlasov-Poisson system, existence of global smooth solutions 
has been proved (see \cite{Pf}); 
in the present case, due to the non-linearity of the Monge-Amp\`ere equation,
we were only able to obtain a result for finite time.

\subsection{The periodic Vlasov-Monge-Amp\`ere system}
\subsubsection*{Polar factorization of maps in a periodic domain}
The polar decomposition Theorem has been generalized in \cite{Mc2} to general Riemannian manifolds,
while the particular case of the flat torus $\Td=\Rd/\mathbb Z^d$ had been addressed in \cite{Co}. 
\begin{defi}
We say that a mapping $Y: \Rd \rightarrow \Rd$ is 
$\mathbb Z^d$ additive if the mapping $x\rightarrow Y(x)-x$ is $\mathbb Z^d$ periodic. 
The set of all measurable $\mathbb Z^d$ additive mappings is denoted by ${\cal P}$. 
For each $x\in \Rd$ we call $\hat x$ 
the class of $x$ in $\Rd/\mathbb Z^d$, and for any $X\in {\cal P}$, $\hat X$ the mapping 
of $\Td$ into itself defined by
$$\forall x \in \Rd, \hat X(\hat x)=\hat{X(x)}.$$ 
\end{defi}
We may say if no confusion is possible additive instead of $\mathbb Z^d$ additive. 
Then the  following theorem can be deduced from the results of \cite{Co} and \cite{Mc2}:

\begin{theo}\label{1McCann}
Let $X: \Rd \rightarrow \Rd$ be additive and 
assume that $\rho_X= X\#dx$ has a density  in  $L^1([0,1]^d)$.
Then there exists a unique pair $(\nabla\Phi_X, \pi_X)$ such that
$$X=\nabla\Phi_X\circ \pi_X$$ 
where
$\Phi_X$ is a convex function and
$\Phi_X(x)-|x|^2/2$ is $\Zd$ periodic,
$\pi_X: \Rd \rightarrow \Rd$ is additive and
$\hat \pi_X$ is Lebesgue measure preserving in $\Td$.  
Moreover we have 
$$\|X- \pi_X\|_{L^2([0,1]^d)}=\|\hat X- \hat \pi_X\|_{L^2(\Td)}
$$ 
and, $\Psi_X$ denoting the Legendre transform of $\Phi_X$, we have
$$\pi_X=\nabla\Psi_X\circ X.$$
\end{theo}

\noindent
Remark : The pair $(\Phi_X, \Psi_X)$ is uniquely defined by the density 
$\rho_X=X\#dx$.
 
\medskip
\noindent
Notice that the periodicity of $\Phi_X(x)-|x|^2/2$
implies that $\nabla\Phi_X$ and $\nabla\Psi_X$ are $\mathbb Z^d$ additive, and that 
$\Psi_X -|x|^2/2$ is also $\mathbb Z^d$ periodic.
As in the previous case, we introduce the following notation:
\begin{defi}\label{1MaTdrho}
Let $\rho$ be a probability  measure on $\Td$, with density in $L^1(\Td)$. We denote $\Phi[\rho]$ 
(resp. $\Psi[\rho]$) the unique up to a constant convex function such that
 \beq
&& \Phi[\rho]-|x|^2/2 \mbox{  is } \Zd \mbox{  periodic },\\ 
&&\forall f \in C^0(\Td),\;\int_{\Td} f(\hat{\nabla\Phi}[\rho](x))dx=\int fd\rho
\enq
(resp. its Legendre fenchel transform).
\end{defi}
$\Psi[\rho]$ will thus be a generalized 
solution of the Monge-Amp\`ere equation \\
$\det D^2\Psi=\rho.$
\\Next the results of Caffarelli (\cite{Ca1}, \cite{Ca2}, \cite{Ca3})
 on the regularity of solutions 
of the  Monge-Amp\`ere equation yield the following theorem:

\begin{theo}\label{1regper}
Let $\rho>0$ be a $ C^{\alpha}(\Td)$ 
probability density on $\Td$, for some $\alpha \in ]0,1[$.
\\
Then $\Psi=\Psi[\rho]$ (see Definition \ref{1MaTdrho}) is a classical solution of 
\begin{eqnarray*} 
\det D^2\Psi=\rho\label{1mongeper}
\end{eqnarray*}
and satisfies:
\begin{eqnarray*} 
&&\|\nabla\Psi(x)-x\|_{\Linf}\leq C(d)=\sqrt d /2\\
&&\|D^2\Psi\|_{C^{\alpha}}\leq K(m,M,\|\rho\|_{C^{\alpha}})
\end{eqnarray*}
where $m=\inf\rho$ and $M=\sup\rho$.
\end{theo}
This theorem is an adaptation of the regularity results stated above, whose complete proof 
is given in appendix.
\subsubsection*{The periodic Vlasov-Monge-Amp\`ere system}
We now seek $f: (t,x,\xi) \in (\Td\times \Rd \times [0,T]) \rightarrow f(t,x,\xi)\in \R^+$, for some $T>0$, 
solution of the initial value problem for the periodic Vlasov-Monge-Amp\`ere $(VMA_p)$ system  
\beq
&&\frac{\partial f}{\partial t}+\nabla_x\cdot\left(\xi f\right)
+\frac{1}{\epsilon^2}\nabla_{\xi}\cdot\left((\nabla\Psi[\rho](x) - x)f\right)=0\label{1principaleper}\\
&&f(0,\cdot,\cdot)=f^0\label{1prinitper},
\enq
for a given $f^0$ satisfying  the compatibility condition
\beq
\int  f^0(x,\xi)dxd\xi= 1.\label{1jaugeper}
\enq
The macroscopic density $\rho$ is still related to $f$ by equation (\ref{1defrho}), and $\Psi[\rho]$ is as in Definition \ref{1MaTdrho}.

\subsection{Existence of local strong solutions} 

We mention first that the proof of existence of global weak solutions adapts
with minor changes to the periodic case, and that the obtain for the periodic $(VMA_p)$ system the same result as Theorem \ref{1globalweak}. 

Our result in this section is the following:
\begin{theo}\label{1strong}
Let  $f_0 \in W^{1,\infty}(\Td\times\Rd),$  be such that:
\beq
&&\exists C_0 >0:\; f_0 \equiv 0\text{ for } |\xi|\geq C_0,\label{1pasvide1}\\
&&\exists m >0:\;\rho_0 (x)=\int_{\Rd}f_0(x,\xi)d\xi\geq m\;\forall x\in \Td\label{1pasvide2},
\enq
then there exists $T>0$ and a solution $f$ to the $VMA_p$ system
(\ref{1principaleper},\ref{1prinitper}), in the space 
$W^{1,\infty}([0,T]\times\Td \times \Rd)$.
\end{theo}

\bigskip
\noindent
{\bf Proof of Theorem \ref{1strong}:}
First we deduce from Theorem \ref{1regper}:
\begin{cor}\label{1regularite}
Let $\rho,\Psi=\Psi[\rho]$ be as in Theorem \ref{1regper}. Then, we have
$$\|D^2\Psi\|_{\Linf(\Td)}\leq C(m,M,\|\nabla\rho\|_{\Linf(\Td)}),$$
and we can define 
$$K(m,M,l)=\sup\{\|D^2\Psi[\rho] \|_{L^{\infty}(\Td)};
\;\;\|\nabla\rho\|_{\Linf(\Td)}\leq l,\;\;\; m\leq\rho\leq M\}<\infty.$$ 
\end{cor}
We see that in order to use Theorem \ref{1regper} we need $\rho$ to be bounded
away from 0. 
In the following lemma, we show that under suitable assumptions on
 the initial data, it is possible to enforce locally in time
 the condition $0<m\leq\rho\leq M.$

\begin{lemme}\label{1densite}
Let $f\in \Linf([0,T]\times\Td\times\Rd)$ satisfy
\begin{eqnarray*}
&&\frac{\partial f}{\partial t}+\nabla_x\cdot\left(\xi f\right)
+\nabla_{\xi}\cdot\left(E(t,x) f\right)=0\\
&&f(0,.,.)=f^0,
\end{eqnarray*}
with $E\in L^1([0,T],BV(\Td))$ and 
\be
\|E\|_{\Linf([0,T]\times\Td)}\leq F,
\en
let the initial condition $f_0$ be such that
$$a(x,\xi)\leq f(0,x,\xi)\leq b(x,\xi),$$
with $\rho_a(x)=\int a(x,\xi)d\xi\geq m>0$ and 
$\rho_b(x)=\int b(x,\xi)d\xi\leq M<\infty$ and $a,b$ satisfying 
\be
|\nabla_{x,\xi}(a,b)|\leq \frac{c}{1+|\xi|^{d+2}}.
\en
Then there exists a constant $R>0$ depending on $m,M,c,F$, such that
\be
(\rho_a(x)-Rt)\leq\rho(t,x)\leq(\rho_b(x)+Rt).
\en
\end{lemme}
The proof of the lemma is given in appendix.

\subsubsection{Construction of approximate solutions}

Let us consider $(t,x)\rightarrow E(t,x)$ a smooth vector-field on $\Td$, and write
\be
T_E(f)=\frac{\partial f}{\partial t}+\nabla_x\cdot\left(\xi f\right)
+\nabla_{\xi}\cdot\left( E\,f\right).
\en
If $f$ satisfies $T_E(f)=0$, we have 
\be
&&T_E\nabla_xf=-(\nabla_xE)\cdot\nabla_{\xi}f\label{1Tx}\\
&&T_E\nabla_{\xi}f=-\nabla_xf\label{1Tv}\\
&&T_E\partial_tf=-\partial_tE\cdot\nabla_{\xi}f,
\en
and therefore 
\beq\label{1Gronwall}
\frac{d}{dt}\|\nabla_{x,\xi}f\|_{\Linf}\leq \|\nabla_{x,\xi}f\|_{\Linf}(1+\|\nabla_x E\|_{\Linf})
\enq
which implies 
\be
\|\nabla_{x,\xi}f(t)\|_{\Linf}\leq \|\nabla_{x,\xi}f(t=0)\|_{\Linf}\exp \left(\int_{0}^{t}(1+\|\nabla_x E(s)\|_{\Linf})ds\right).
\en

\medskip
\noindent
Now let $f_0$ be given as in Theorem \ref{1strong}, satisfying (\ref{1pasvide1},\ref{1pasvide2}).  
Thanks to Lemma \ref{1densite} it is possible to find $t_1,m,M$ such that for any $f$ satisfying 
\be
&&T_E(f)=0\\
&& f(t=0)=f_0,
\en
for any field $E\in L^1([0,t_1],BV(\Td))$ satisfying
$\|E\|_{\Linf([0,t_1]\times \Td)}\leq \sqrt{d}/(2\epsilon^2)$, 
we have
\beq
&&m\leq\rho(t,\cdot)\leq M,\;\;\forall t\in [0,t_1]\label{1cond1}\\
&&|\xi|_{max} \leq C_1=10 C_0\label{1cond2},
\enq with $f$ supported in $\{|\xi|\leq |\xi|_{max}\}$ and with $C_0$ as in 
Theorem \ref{1strong}, so that we have for $0\leq t\leq t_1$: 
\beq\label{1nablarhof}
\|\nabla\rho\|_{\Linf}\leq \omega_d C_1^d\|\nabla_x f\|_{\Linf},
\enq 
$\omega_d$ being the volume of the unit ball of $\Rd$.
Then we construct a family of approximate solutions $(f_h,\Psi_h)$
to (\ref{1principaleper}),
in the same spirit as we did for weak solutions,
by solving
\be
&& \frac{\partial f_h}{\partial t}+\xi\cdot\nabla_xf_h
+\frac{\nabla\Psi_h(x)-x}{\epsilon^2}\cdot\nabla_{\xi} f_h=0\\
&&f_h(t=0)=f_0\\
&&\Psi_h(t)=\Psi[\rho(t=nh)] \mbox{ for } t\in [nh,(n+1)h[.
\en
Note that we have neither mollified the term $\nabla\Psi$ nor the initial condition and that $\|\nabla\Psi_h-x\|_{\Linf}\leq C(d)=\sqrt{d}/2$.
Now choose $l=10\|\nabla_{x,\xi}f_0\|_{\Linf}\omega_d C_1^d$, if for some $t=nh\leq t_1-h$ we have 
\be
\|\nabla_{x,\xi}f^h(t=nh)\|_{\Linf} \leq \frac{l}{\omega_d C_1^d}
\en
this implies, thanks to (\ref{1nablarhof}), that
\be
\|\nabla_x\rho^h(t=nh)\|_{\Linf}\leq l,
\en
and conditions (\ref{1cond1},\ref{1cond2}) are satisfied because $t\leq t_1$. 
Then if we denote $K=K(m,M,l)$ as in Corollary \ref{1regularite}, we have for $nh\leq t< nh+h,$
\be
\frac{d}{dt}\|\nabla_{x,\xi}f^h\|_{\Linf}\leq (K+1)\|\nabla_{x,\xi}f^h\|_{\Linf},
\en
and then  
\be
\|\nabla_{x,\xi}f^h(t=nh+h)\|_{\Linf} \leq \|\nabla_{x,\xi}f^h(t=nh+h)\|_{\Linf}\exp{(K+1)h}.
\en
So if we define $T$ as 
\be
T=\min \{t_1,t_2\},
\en
with $\exp((K+1)t_2)=10$,
we have for $0\leq t \leq T$,
\be
&&\|\nabla_{x,\xi}f^h\|_{\Linf}\leq 10\|\nabla_{x,\xi}f_0\|_{\Linf}\\
&&\|\nabla\rho^h\|_{\Linf}\leq l\\
&&m\leq\rho\leq M\\
&&\|D^2\Psi^h\|_{\Linf}\leq K.
\en

 Thus we can extract a subsequence converging to a strong solution of (\ref{1principaleper},\ref{1prinitper}). Then we argue as in section 2 to show that all terms converge to the correct limit. This ends the proof of Theorem \ref{1strong}.

$\hfill\Box$


\section{Asymptotic analysis}
\subsection{Convergence to the Euler equation}\label{1asymptotic}
In this section we justify that the Vlasov-Monge-Amp\`ere system describes approximate geodesics on volume preserving transformations: indeed we will show that 
weak solutions of this system converge to a solution 
of the incompressible Euler equations $(E)$ as the parameter $\epsilon$ 
goes to 0, at least for well prepared initial data. 
We restrict ourselves to the space periodic case, the macroscopic density $\rho$ is still defined by 
(\ref{1defrho}) and the convex potentials $\Phi[\rho],\Psi[\rho]$ are still as in Definition \ref{1MaTdrho}.


For sake of simplicity, we slightly modify our notations and introduce the following rescaled potentials
\be
&&\tilde\varphi[\rho]=\frac{|x|^2/2-\Psi[\rho]}{\epsilon},\\
&&\varphi[\rho]=\frac{\Phi[\rho]-|x|^2/2}{\epsilon},
\en
so that 
\be
\nabla\varphi[\rho]=\nabla\tilde\varphi[\rho]\circ \nabla\Phi[\rho],
\en
and the $VMA_p$ system takes the following form: 
\begin{eqnarray}
&&\frac{\partial f}{\partial t}+\xi\cdot\nabla_x f
-\frac{\nabla\tilde\varphi[\rho]}{\epsilon}\cdot\nabla_{\xi}f=0
\label{1principale2}\\
&&f(0,\cdot,\cdot)=f_0\label{1initiale}.
\end{eqnarray}
 The energy is given by
\begin{eqnarray}\label{1energyvma}
E(t)&&=\frac{1}{2}\int f(t,x,\xi)|\xi|^2 dxd\xi + 
\frac{1}{2}\int|\nabla\varphi|^2 dx\\
&&=\frac{1}{2}\int f(t,x,\xi)|\xi|^2 dxd\xi + 
\frac{1}{2}\int\rho|\nabla\tilde\varphi|^2 dx\nonumber.
\end{eqnarray}
It has been shown in section \ref{1energysection} that the energy is conserved. The Euler equations for incompressible fluids $(E)$ reads:
\beq
&&\dt v + v\cdot\nabla v =-\nabla p\label{1euler1}\\
&& \nabla\cdot v = 0\label{1euler2}.
\enq
We shall here consider a smooth solution of $E$ and a weak solution of
$VMA_p$, with `well prepared initial data',
meaning that the initial data of both systems are close a time 0. 
Then we will show that as time evolves, both solutions stay close to each other.
\begin{theo}
\label{1neutre}
Let $f$ be a weak solution of (\ref{1principale2}, \ref{1initiale}) with finite energy, let 
$(t,x)\rightarrow v(t,x)$ be a smooth $C^2([0,T]\times\Td)$ solution of (\ref{1euler1},\ref{1euler2}) for $t\in[0,T],$ and $p(t,x)$ the corresponding pressure, let  
\begin{eqnarray*}
H_\epsilon(t)=\frac{1}{2}\int f(t,x,\xi)|\xi-v(t,x)|^2 dx d\xi + \frac{1}{2}\int|\nabla\varphi(t,x)|^2 dx,
\end{eqnarray*}
then 
\begin{eqnarray*}
H_{\epsilon}(t)\leq  C\exp(Ct)(H_{\epsilon}(0) + \epsilon^2),\;\forall t \in [0,T].
\end{eqnarray*}
$C$  depends only on 
$T, \sup_{0\leq s\leq T} \left\{ \|v(s,.), p(s,.), \partial_t p(s,.), \nabla p(s,.)\|_{W^{1,\infty}(\Td)}\right\}$. 

\end{theo}

\medskip
\noindent
Remark 1: This estimate is enough to compare the weak solutions $f$ to
the $VMA_p$ system (for well prepared initial data) and the smooth solutions $v$ of the Euler equations.
For instance, $\int f(t=0,x,\xi)d\xi\equiv 1$
implies \mbox{$\varphi(t=0,x)\equiv 0$} and therefore,
$$\int|\xi-v(t=0,x)|^2 f(t=0,x,\xi) dxd\xi\le C_0\epsilon^2
$$
implies
$$
\sup_{t\in [0,T]}
\int|\xi-v(t,x)|^2 f(t,x,\xi)dxd\xi\le C_T\epsilon^2,
$$
where $C_T$ depends only on $C_0$, $T$ and $v$.
\\
Remark 2: We see that we consider nearly monokinetic initial data for
the $VMA_p$ system. 
\subsubsection*{Proof of Theorem \ref{1neutre}} 
We shall show that
\begin{eqnarray}\label{1dtH} 
\frac{d}{dt}H_{\epsilon}   =
&&-\int f(t,x,\xi)(\xi-v)\nabla v(\xi-v)\nonumber\\
&&+\int f(t,x,\xi)\frac{1}{\epsilon}v\cdot\nabla\tilde\varphi\nonumber\\
&&-\int f(t,x,\xi)(v-\xi)\cdot\nabla p,
\end{eqnarray}
where we will use the notation 
\be
u  \ \nabla v \  w= \sum_{i,j=1}^d u^i \partial_{i}v^j w^j.
\en
The proof of this identity is postponed to the end of the section.

Now we look at all terms of the right hand side. All the constants that we denote by $C$ are controlled as in Theorem \ref{1neutre}. 
We set 
\be
&&T_1=-\int f(t,x,\xi)(\xi-v)\nabla v(\xi-v),\\
&&T_2=\int f(t,x,\xi)\frac{1}{\epsilon}v\cdot\nabla\tilde\varphi,\\
&&T_3=-\int f(t,x,\xi)(v-\xi)\cdot\nabla p.
\en
First we have $T_1\leq C H_{\epsilon}.$
For $T_2$ we have 
\begin{eqnarray*}
T_2=\frac{1}{\epsilon} \int \rho v \cdot\nabla\tilde\varphi =&&
 \frac{1}{\epsilon}\int v(\nabla\Phi[\rho])
\cdot\nabla\tilde\varphi(\nabla\Phi[\rho])\\ 
=&& \frac{1}{\epsilon}\int v(x+\epsilon\nabla\varphi)\cdot\nabla\varphi\\ 
=&& \frac{1}{\epsilon}\int v\cdot\nabla\varphi + (v(x+\epsilon\nabla\varphi)
-v(x))\cdot\nabla\varphi\\
\leq &&0 + C \int\left|\nabla\varphi\right|^2 \leq C H_{\epsilon},  
\end{eqnarray*}  
we have used that $v$ is divergence-free thus its integral against any gradient is zero.
Next we have the following lemma: 
\begin{lemme}
\label{1rofi}Let  $G:\Td\rightarrow \R$ be Lipschitz continuous such that $\displaystyle \int_{\Td} G =0$, then for all $R>0$, one has
$$|\int \rho G| \leq   \frac{1}{2}\|\nabla G\|_{\Linf}(\frac{1}{R}\epsilon^2 +  R H_{\epsilon}).$$  
\end{lemme}
{\bf Proof:} We just write a Taylor expansion of $G$:   
\begin{eqnarray*}
&&\left|\int (\rho-1) G\right|=\left|\int(G(x+\epsilon\nabla\varphi)-G(x)\right|\\
\leq &&\epsilon\|\nabla G\|_{\Linf} \|\nabla\varphi\|_{L^1}\leq \epsilon\|\nabla G\|_{\Linf} H_{\epsilon}^{1/2}\leq \frac{1}{2}\|\nabla G\|_{\Linf} (\frac{1}{R}\epsilon^2 + R H_{\epsilon}).
\end{eqnarray*}
$\hfill$ $\Box$

\bigskip
\noindent
Again, since $v$ is divergence-free, $\int v\cdot\nabla p=0,$ thus from Lemma \ref{1rofi} we have
\be
-\int\rho v\cdot\nabla p\leq C (\epsilon^2 + H_{\epsilon})\nonumber.
\en 
We remind that $$\dt \rho(t,x) = -\nabla_x\cdot\int f(t,x,\xi)\xi d\xi.$$
Since it costs no generality to suppose that for all $t\in [0,T]$, 
$\int p(t,x)dx \equiv 0$, we obtain that
\be
\int f(t,x,\xi)\xi\cdot\nabla p
&=&\int \frac{\partial\rho}{\partial t} p\\
&=&\frac{d}{dt}\int\rho p - \int \rho\frac{\partial p}{\partial t}\nonumber\\
&\leq & C(\epsilon^2 + H_{\epsilon})-\frac{dQ}{dt} \nonumber\\
\en
again using Lemma \ref{1rofi}, where $\;\displaystyle Q(t)=-\int \rho p.$ Thus 
\be 
T_3\leq C(H_{\epsilon}+\epsilon^2) - \frac{dQ}{dt} 
\en
and we have the following inequality:
\beq
\frac{d}{dt}(H_{\epsilon}+Q)\leq CH_{\epsilon}+ O(\epsilon^2).\label{1qwe}
\enq
Moreover, using Lemma \ref{1rofi},
\beq
|Q(t)|\leq C\epsilon^2 + H_{\epsilon}(t)/2,\label{1Q}
\enq
thus
\beq
H_{\epsilon}+Q\geq H_{\epsilon}/2- C\epsilon^2,
\enq
and we can transform (\ref{1qwe}) in 
\begin{eqnarray}
\frac{d}{dt}(H_{\epsilon}+Q)\leq C(H_{\epsilon}+Q)+ C\epsilon^2.
\end{eqnarray}
Gronwall's lemma then yields
$$H_{\epsilon}(t)+Q(t)\leq (H_{\epsilon}(0)+Q(0)+Ct\epsilon^2)\exp(Ct).$$
Using again  (\ref{1Q}) we obtain 
\begin{eqnarray}
H_{\epsilon}(t)\leq C(H_{\epsilon}(0)+\epsilon^2)\exp(Ct),
\end{eqnarray}
which achieves the proof of Theorem \ref{1neutre}.
 
$\hfill$ $\Box$

\bigskip
\noindent
{\bf Proof of identity (\ref{1dtH}):}

\noindent 
We first notice that, 
for all $g\in C^1(\mathbb R\times\Td)$ , we have:
$$\frac{d}{dt}\int\rho(t,x) g(t,x) dx=\int\int f(t,x,\xi)(\partial_t g(t,x) + \xi\cdot\nabla g(t,x)) d\xi dx.$$  We also use the conservation of energy defined by (\ref{1energyvma}). Then we get
\begin{eqnarray*}
\frac{d}{dt}H_{\epsilon}=&&\frac{d}{dt}\frac{1}{2}\int f(t,x,\xi)(|v|^2- 2\xi\cdot v)dx d\xi\\
=&&\int f(t,x,\xi)(\dt v\cdot v - \dt v \cdot\xi)
- \frac{1}{2}\int \nabla_x\cdot(f(t,x,\xi)\xi)(|v|^2- 2\xi\cdot v) \\ 
&&  + \frac{1}{2}\int \nabla_\xi\cdot(\frac{1}{\epsilon}\nabla\tilde\varphi f(t,x,\xi))(|v|^2- 2\xi\cdot v).
\end{eqnarray*}
Integrating by part, we get
\begin{eqnarray*}
\frac{d}{dt}H_{\epsilon}=&&\int f(t,x,\xi)(\dt v\cdot v - \dt v \cdot\xi)+ \int f(t,x,\xi)\xi \nabla v (v-\xi)\\ 
&&+\int f(t,x,\xi)\frac{1}{\epsilon}\nabla\tilde\varphi\cdot v.
\end{eqnarray*}
The first two terms can be rewritten as
\begin{eqnarray*}
&&\int f(t,x,\xi)(\dt v\cdot v - \dt v \cdot\xi)+ 
\int f(t,x,\xi)\xi \nabla v (v-\xi)\\ 
=&&-\int f(t,x,\xi)(v-\xi) \nabla v (v-\xi)
+\int f(t,x,\xi)\dt v\cdot(v - \xi)\\
&& + \int f(t,x,\xi)v \nabla v (v-\xi)\\
=&&-\int f(t,x,\xi)(v-\xi) \nabla v (v-\xi) 
+ \int f(t,x,\xi)(v-\xi)\cdot(\dt v + v\cdot\nabla v),
\end{eqnarray*} and finally using equation (\ref{1euler1}) we conclude.

$\hfill$ $\Box$

\subsection{Comparison with the Euler-Poisson system}\label{1section-Euler-Poisson}
Here we show that, as mentioned in the introduction, 
the Euler-Poisson ($EP$) system
is a more accurate 
approximation to the Vlasov Monge-Amp\`ere
system than the Euler equations, as $\epsilon$ goes to zero.

\bigskip
\noindent
{\bf The $EP$ system}
Let us recall
that the (pressureless) Euler-Poisson system
describes the motion of a continuum of electrons 
on a neutralizing background of ions through electrostatic interaction.
Let $\bar v$ and $\bar\rho$ be the velocity and density of electrons. 
Let $\bar\varphi$ be the (rescaled) electric potential.
Under proper scaling, these functions of $x\in \Rd$ and $t>0$ satisfy the Euler-Poisson 
system:
\beq
&&\dt \bar v + \bar v\cdot\nabla \bar v =-\frac{1}{\epsilon}\nabla\bar\varphi
 \label{1eulerpoisson1}\\
&&\dt \bar\rho + \nabla\cdot(\bar\rho \bar v)=0\label{1eulerpoisson2}\\
&& 1-\epsilon\Delta\bar\varphi=\bar\rho\label{1eulerpoisson3}.
\enq
The so-called 'quasi-neutral' limit $\epsilon\rightarrow 0$ of similar systems 
has been studied for example in \cite{Gr}  and \cite{CG}, and convergence
results have been established using pseudo-differentials energy estimates.
For well-prepared initial data, solutions of $EP$ are expected to behave
as solutions of Euler incompressible equations.
This fact is proved by the second author in his PhD thesis (\cite{These}, Chap 2), see also \cite{L1}. 
We give here the complete result that we will use herafter.
We will denote by $\bar v\ep$ (resp. $f\ep)$ the solutions
of the $EP$ (resp. $VMA_p$) system with parameter $\epsilon$.
\begin{theo}\label{1cvep}
Let $v$ be a solution of (\ref{1euler1},\ref{1euler2})
on $[0,T]\times \Td$,
with initial data $v_0$, and satisfying $v\in \Linf([0,T], H^s(\Td))$ for some $s\geq s_0(d)$. There for some $s'>0$, $s'< s$, if
$(\bar v_0\ep, \bar\rho_0\ep)$ is such that the sequences
\be
\frac{\bar v_0\ep-v_0}{\epsilon}, \hspace{1cm}  \frac{\bar\rho_0\ep-1}{\epsilon^2} 
\en
are bounded in $H^{s'}(\Td)$,  then there exists $T\epu>0$ with $\liminf_{\epsilon\rightarrow 0}T_{\epsilon}\geq T$ and a sequence $(\bar v\ep,\bar\rho\ep)$ of solutions to the $EP$ system on $[0,T_{\epsilon}[$
with initial data $(\bar v_0\ep, \bar\rho_0\ep)$,  belonging to $\Linf([0,T\epu], H^{s'}(\Td))$.  Moreover, for $\epsilon$ 
small enough, the sequences
\be
\frac{\bar v\ep-v}{\epsilon}, \hspace{1cm} \frac{\bar\rho\ep-1}{\epsilon^2}
\en
are bounded in $\Linf([0,T], H^{s'}(\Td))$. Finally,  $s'$ goes to $+\infty$ as $s$ goes to $+\infty$. 
\end{theo}

\subsubsection*{Assumptions}
Here we consider $v$ a solution to $E$ (\ref{1euler1}, \ref{1euler2}) 
with initial data $v_0$, a sequence $f\ep$ of  solutions of $VMA_p$
(\ref{1principale2},\ref{1initiale}) with initial data $f\ep_0$,  and a sequence $(\bar v\ep, \bar\rho\ep)$ solutions of
$EP$ (\ref{1eulerpoisson1}, \ref{1eulerpoisson2}, \ref{1eulerpoisson3}) with initial data $(\bar v \ep_0, \bar\rho\ep_0)$.
We still define $H\epu$ as in Theorem \ref{1neutre}:
\be
H\epu(t)=\frac{1}{2}\int f\ep(t,x,\xi)|\xi- v(t,x)|^2 dx d\xi 
+ \frac{1}{2}\int|\nabla\varphi\ep|^2 dx. 
\en
We introduce the following assumptions:

\begin{enumerate}
\item[{\bf H0}] $v$ solution of $E$ satisfies, for some $C_0>0$,  $\ds\|v\|_{\Linf([0,T] H^s(\Td))} \leq C_0$, and $s$ is large enough so that $s'$ in Theorem \ref{1cvep} satisfies $\ds s \geq  s'>[\frac{d}{2}] + 2$.

\item[{\bf H1}] The sequence $(\bar v_0\ep , \bar\rho_0\ep)$ of initial data  of $EP$
is such that, for some $C_1>0$, 
\be
\sup_{\epsilon > 0}\left\{\frac{1}{\epsilon} \|\bar v\ep_0-v\|_{H^{s'}(\Td)}, \  
\frac{1}{\epsilon^2}\|\bar\rho\ep-1\|_{H^{s'}(\Td)}\right\} \leq C_1. 
\en

\item[{\bf H2}] The sequence $f\ep_0$ satisfies $H\epu(0)\leq C_2\epsilon^2$ for some $C_2>0$.

\end{enumerate}

\bigskip
\noindent
{\bf H0}, {\bf H1}, {\bf H2} imply that
\begin{enumerate} 
\item There exists $\tilde C_0$ such that 
\beq\label{1boundei}
\|v\|_{\Linf([0,T],W^{2,\infty}(\Td))}\leq \tilde C_0.
\enq 
\item From Theorem \ref{1neutre}, there exists  $\tilde C_1$ such that
\beq\label{1boundvma}
H\epu(t)\leq  \tilde C_1\epsilon^2 \mbox{ for } t\in[0,T].
\enq
\item From Theorem \ref{1vmavp} and Sobolev imbeddings, there exists  $\tilde C_2$ such that
\beq\label{1boundep}
\sup_{\epsilon < \epsilon_0} \left\{ \left\|\frac{\bar v\ep-v}{\epsilon}, \hspace{.5cm} \frac{\bar\rho\ep-1}{\epsilon^2} \right\|_{\Linf([0,T],W^{2,\infty}(\Td))}\right\} \leq \tilde C_2. 
\enq
\end{enumerate}

\bigskip
\noindent
We are now ready to prove the following result:

\begin{theo}
\label{1vmavp}
Let $f\ep_0,\bar v\ep_0,\bar\rho\ep_0, v, T$ be as above, satisfying assumptions {\bf H0}, {\bf H1}, {\bf H2}.
Define
\begin{eqnarray*}
G_\epsilon(t)=\frac{1}{2}\int f\ep(t,x,\xi)|\xi-\bar v\ep(x)|^2 dx d\xi
 + \frac{1}{2}\int|\nabla\varphi\ep-\nabla\bar\varphi\ep|^2 dx.
\end{eqnarray*}
Then there exists $C>0$  such that 
\begin{eqnarray*}
G_{\epsilon}(t)\leq C\exp(Ct)(G_{\epsilon}(0) + \epsilon^3),\;\forall t \in [0,T]
\end{eqnarray*}
where C depends on $s', C_0, C_1, C_2, T$.

\end{theo}
\bigskip
\noindent
Remark: the theorem shows that the distance between solutions of the $(EP)$ system 
and the $VMA_p$ system measured with $G\epu$ is like $O(\epsilon^3)$ 
whereas Theorem \ref{1neutre} showed
that the distance between the solution of
the Euler equation and the $VMA_p$ system was like $O(\epsilon^2)$.
Note also that $G\epu$ and $H\epu$ can both be interpreted as the square of a distance.

\bigskip
\noindent
{\bf Proof of Theorem \ref{1vmavp}:}
For notational simplicity, we drop most $\epsilon$'s. 
Proceeding as in (\ref{1dtH}) and noticing that:
\be
\Dt\int_{\Td}|\nabla\bar\varphi|^2=\frac{1}{\epsilon}\int_{\Td}\bar \rho \bar v
\cdot\nabla\bar\varphi
\en
we obtain the following identity:
\begin{eqnarray}\label{1dtH2} 
\frac{d}{dt}G_{\epsilon}   =
&&-\int f(t,x,\xi)(\xi-\bar v)\nabla \bar v(\xi-\bar v)\nonumber\\
&&+\int f(t,x,\xi)\frac{1}{\epsilon}\bar v\cdot\nabla\tilde\varphi
-\int f(t,x,\xi)\frac{1}{\epsilon}\bar v \cdot\nabla \bar\varphi       \nonumber\\
&&+\int f(t,x,\xi)\frac{1}{\epsilon}\xi\cdot\nabla \bar\varphi
+\int \frac{1}{\epsilon}\bar\rho\bar v\cdot\nabla\bar\varphi\nonumber\\
&&- \Dt \int \nabla\bar\varphi\cdot\nabla\varphi.
\end{eqnarray}
Then we notice
\be
\int f(t,x,\xi)\frac{1}{\epsilon}\xi\cdot\nabla \bar\varphi
=\Dt\left(\int \frac{1}{\epsilon} \rho\bar\varphi\right)-\frac{1}{\epsilon}\int\rho\dt\bar\varphi.
\en
Next, we have the following lemma:
\begin{lemme}\label{1Taylor}
Define for any $\theta\in C^2(\Td)$
\be
&&<\nabla\theta>(x)=\int_{0}^{1}\nabla\theta(x+s\epsilon\nabla\varphi(x))ds,\\
&&<\nabla^2 \theta>(x)=\int_{0}^{1}(1-s)\nabla^2 \theta(x+s\epsilon\nabla\varphi(x))ds.
\label{1def<>}
\en
Then 
\be
\int \rho\theta &=&\int\theta+\epsilon\int <\nabla\theta>\cdot\nabla\varphi \\
&=& \int\theta+\epsilon\int \nabla\theta\cdot\nabla\varphi + \epsilon^2\int <\nabla^2 \theta>\nabla\varphi\nabla\varphi.
\en
\end{lemme}
{\bf Proof:} The proof just uses the Taylor expansion  and the identity \\
$\int \rho\theta=\int\theta(x+\epsilon \nabla\varphi)$.

$\hfill \Box$

\bigskip
\noindent
Using Lemma \ref{1Taylor}, we get
\be
&&\frac{1}{\epsilon}\int \rho \dt\bar\varphi \\
=&&\frac{1}{\epsilon}\int\dt\bar\varphi+\int\dt\nabla\bar\varphi\cdot\nabla\varphi
+\epsilon\int<\dt \nabla^2\bar\varphi>\nabla\varphi\nabla\varphi.
\en
We claim that, under our assumptions, we have 
$$\|\dt \nabla^2\bar\varphi\|_{\Linf([0,T']\times\Td)}\leq C.$$

\noindent
Proof: from (\ref{1eulerpoisson2}), we have 
\be
\dt\bar\rho=-\bar\rho\nabla\cdot \bar v - \bar v\cdot \nabla\bar\rho.
\en
Using (\ref{1boundep}), we obtain that $\ds\|\dt\bar\rho\|_{H^{s'-1}}\leq C \epsilon$. Since $H^{s'}(\Td)$ is continuously embedded in
$W^{2,\infty}(\Td)$,  $H^{s'-1}(\Td)$ is continuously embedded in $\Linf(\Td)$.\\ 
Then, using (\ref{1eulerpoisson3})
and classical elliptic regularity, we have  
$$\epsilon\|\dt \nabla^2\bar\varphi\|_{H^{s'-1}}\leq C \|\dt\bar\rho\|_{H^{s'-1}},$$
and the desired result follows.

$\hfill\Box$

\bigskip
\noindent
This implies, using (\ref{1boundvma}), that  
\be
\left | \epsilon\int<\dt \nabla^2\bar\varphi>\nabla\varphi\nabla\varphi \right |\leq C\epsilon^3.
\en
Next, 
\be
&&\int\dt\nabla\bar\varphi\cdot\nabla\varphi=-\int \dt\Delta\bar\varphi\varphi\\
&&=\frac{1}{\epsilon}\int \dt\bar\rho\varphi=\frac{1}{\epsilon}\int\bar\rho\bar v \cdot\nabla\varphi.
\en
Using again Lemma \ref{1Taylor}, we get
\be
&&\Dt \int \nabla\bar\varphi\cdot\nabla\varphi\\
=&&\frac{1}{\epsilon}\Dt\left(\int \rho \bar\varphi-
\epsilon^2\int<\nabla^2\bar\varphi>\nabla\varphi\nabla\varphi\right)
\en
and for the same reasons we have $\|\nabla^2\bar\varphi\|_{\Linf([0,T]\times\Td))}\leq C\epsilon$. This yields
$$Q(t)=\epsilon\int<\nabla^2\bar\varphi>\nabla\varphi\nabla\varphi=O(\epsilon^4).$$
Moreover, it does not cost to set $\int \bar\varphi\equiv 0$ and deduce
\be
\int f(t,x,\xi)\frac{1}{\epsilon}\xi\cdot\nabla \bar\varphi- \Dt \int \nabla\bar\varphi\cdot\nabla\varphi
=-\frac{1}{\epsilon}\int\bar\rho\bar v \cdot\nabla\varphi+O(\epsilon^3)+\Dt Q.
\en
Thus the remaining terms are
\be
R=\frac{1}{\epsilon}\int \left[\rho\nabla\tilde\varphi -\rho\nabla\bar\varphi
+\bar\rho\nabla\bar\varphi-\bar\rho\nabla\varphi\right]\cdot\bar v.
\en
Calculations that we postpone to the end of the proof show that 
\beq\label{1dtH3}
R\leq&&\int (\nabla\varphi-\nabla\bar\varphi)\nabla\bar v (\nabla\varphi-\nabla\bar\varphi)\nonumber
+ C\int |\nabla\varphi-\nabla\bar\varphi|^2
\\
&&-\demi\int \nabla\cdot\bar v(|\nabla\bar\varphi|^2-2\nabla\varphi\cdot\nabla\bar\varphi ) + C\epsilon^3.
\enq
with $C$ depending on $\ds\|\nabla^2\bar v \|_{\Linf([0,T]\times\Td)}$ and $\ds\epsilon^{-1}\|\nabla^3\bar\varphi\|_{\Linf([0,T]\times\Td)}$, therefore uniformly bounded thanks to (\ref{1boundep}).
Finally we obtain
\be
\Dt G\epu\leq&&-\int f(t,x,\xi)(\xi-\bar v)\nabla \bar v(\xi-\bar v)+(\nabla\varphi-\nabla\bar\varphi)\nabla\bar v (\nabla\varphi-\nabla\bar\varphi)\\
&& -\demi\int (\nabla\cdot\bar v)(|\nabla\bar\varphi|^2-2\nabla\bar\varphi\cdot \nabla\varphi)+ C\int |\nabla\varphi-\nabla\bar\varphi|^2\\
&& + C\epsilon^3  +  \Dt Q\;\;\;
\en
with $\left|Q(t)\right|\leq C\epsilon^4$ for $t\in[0,T]$.
From (\ref{1boundep}) we have $\ds\|\nabla\cdot v\|_{\Linf([0,T]\times\Td)}\leq C\epsilon$ and $\ds\|\nabla\bar\varphi\|_{\Linf([0,T]\times\Td)}\leq C\epsilon$,
whereas (\ref{1boundvma}) yields $\ds\int \left|\nabla\varphi\right|^2 \leq C\epsilon^2$.
Note that we also have
\be
&&-\int f(t,x,\xi)(\xi-\bar v)\nabla \bar v(\xi-\bar v)+(\nabla\varphi-\nabla\bar\varphi)\nabla\bar v (\nabla\varphi-\nabla\bar\varphi)\\
&&+ C\int |\nabla\varphi-\nabla\bar\varphi|^2\leq C G\epu.
\en
We conclude that
\be
\Dt (G\epu-Q)\leq C ((G\epu-Q)+ \epsilon^3),
\en
and the conclusion of Theorem \ref{1vmavp} follows by Gronwall's lemma.

$\hfill \Box$

\bigskip
\noindent
{\bf Proof of identity (\ref{1dtH3}):}
Here we have to compute:
\be
R=\frac{1}{\epsilon}\int \bar v(x+\epsilon\nabla\varphi)\cdot\nabla\varphi - (\bar v \nabla\bar \varphi)(x+\epsilon\nabla\varphi)
+ (1-\epsilon \Delta \bar\varphi)(\bar v \cdot\nabla\bar \varphi - \bar v\cdot\nabla\varphi)
\en
Using Lemma \ref{1Taylor} we have:
\be
R=&&\frac{1}{\epsilon}\int \bar v\cdot\nabla\varphi -\bar v \cdot\nabla\bar \varphi+\bar v \cdot\nabla\bar \varphi - \bar v\cdot\nabla\varphi\\
&&+\int \nabla\bar v \cdot\nabla\varphi \nabla\varphi -\nabla(\bar v \nabla\bar \varphi)\nabla\varphi 
- \bar v \nabla\bar \varphi\Delta\bar \varphi +\bar v \nabla\varphi\Delta\bar \varphi\\
&&+\int(<\nabla\bar v>-\nabla\bar v)\nabla\varphi\nabla\varphi 
- \epsilon<\nabla^2(\bar v \nabla\bar\varphi)>\nabla\varphi\nabla\varphi.
\en
We see that the first line cancels. Then we show that the last line is bounded by $C\epsilon^3$.
\\
This is obvious for the last term since from (\ref{1boundei}, \ref{1boundvma}) we have 
$\ds\|\bar v\|_{W^{2,\infty}} \leq C$, and $\ds \|\nabla \bar\varphi \|_{W^{2,\infty}}\leq C\epsilon$.
\\
Then for the first term we have the following lemma:
\begin{lemme}
\label{1delta}
We define 
$$\Delta=\int( <\nabla\bar v>(x) - \nabla \bar v(x))\nabla\varphi\nabla\varphi dx,$$ 
then one has:
$$\left|\Delta\right|\leq  C\epsilon^{10/3}+ C \int|\nabla\varphi-\nabla\bar\varphi|^2 dx.$$
\end{lemme}
{\bf Proof:}
\noindent
First we show that if $\Theta(R)=\int_{\{\left|\nabla\varphi\right|\geq R\}}\left|\nabla\varphi\right|^2$, 
$$\Theta(R)\leq C \int\left|\nabla\varphi-\nabla\bar\varphi\right|^2 +\frac{C\epsilon^4}{R^2}.$$

\noindent
Proof: $\int \left|\nabla\varphi\right|^2\leq C\epsilon^2$, implies that
$$\mbox{meas}\{\left|\nabla\varphi\right|\geq R\}\leq C(\frac{\epsilon}{R})^2.$$
Since $\left|\nabla\bar\varphi(t,x)\right|\leq \epsilon$ for $(t,x)\in [0,T'x]\times\Td$
\begin{eqnarray*}
&&\Theta(R)\leq\int_{\{\left|\nabla\varphi\right|\geq R\}}\left|\nabla\bar\varphi\right|^2+\int_{\{\left|\nabla\varphi\right|\geq R\}}\left|\nabla\varphi-\nabla\bar\varphi\right|^2\\
&&\leq \frac{C\epsilon^4}{R^2}+\int\left|\nabla\varphi-\nabla\bar\varphi\right|^2.
\end{eqnarray*}

\medskip
\noindent
Then we have
\begin{eqnarray*}
&&\Delta\leq C \Theta(R)+ \int_{\left|\nabla\varphi\right|\leq R}\left|<\nabla\bar v>(x) - \nabla\bar v(x)\right|\nabla\varphi\nabla\varphi\\
\mbox{ with }&& \left|<\nabla\bar v>(x) - \nabla\bar v(x)\right|\leq C \epsilon\left|\nabla\varphi\right|\\
\mbox{thus } &&\Delta\leq C\epsilon\int_{\left|\nabla\varphi\right|\leq R}\left|\nabla\varphi\right|^3 + C\Theta(R)\\
&&\leq C\left(\epsilon R \int\left|\nabla\varphi\right|^2 +\frac{\epsilon^4}{R^2}+\int \left|\nabla\varphi-\nabla\bar\varphi\right|^2\right)\\
&&\leq C\left(\epsilon^3 R +  +\frac{\epsilon^4}{R^2}+\int \left|\nabla\varphi-\nabla\bar\varphi\right|^2\right)
\en
for all R, so for $R=\epsilon^{(1/3)}$ one obtains:
\be
&&\Delta\leq C\epsilon^{10/3}+C\int \left|\nabla\varphi-\nabla\bar\varphi\right|^2.
\end{eqnarray*}
$\hfill \Box$ 

\bigskip
\noindent
Thus we have shown that $R=S+O(\epsilon^3)$,
and $\displaystyle S=\Sigma_{k=1}^{6}T_k$
where each $T_k$ is given by:
\be
&& T_1= \partial_j\bar v_i\partial_j\varphi\partial_i\varphi\\
&& T_2= -\partial_j\bar v_i\partial_j\varphi\partial_i\bar\varphi\\
&& T_3= -\bar v_i\partial_{ij}\bar\varphi\partial_j\varphi\\
&& T_4= \partial_j\bar v_i\partial_j\bar\varphi\partial_i\bar\varphi\\
&& T_5= \bar v_i\partial_{ij}\bar\varphi\partial_j\bar\varphi\\
&& T_6= \bar v_i\partial_{jj}\bar\varphi\partial_i\varphi\\
\en
where we have used Einstein's convention for repeated indices.
First we have
\be
T_5=-\demi\int(\nabla\cdot \bar v) \left|\nabla\bar\varphi\right|^2
\en
\be
T_1+T_2+T_4=\int \partial_j\bar v_i(\partial_j\varphi-\partial_j\bar\varphi)(\partial_i\varphi-\partial_i\bar\varphi)+ T_7
\en
with $T_7=\int \partial_j\bar v_i\partial_j\bar\varphi\partial_i\varphi$.
\be
T_6=-\int \partial_i\bar v_i\partial_{jj}\bar\varphi\varphi + \bar v_i\partial_{ijj}\bar\varphi\varphi
\en
and 
\be
-\int \bar v_i\partial_{ijj}\bar\varphi\varphi=\int \partial_j\bar v_i\partial_{ij}\bar\varphi\varphi
+\bar v_i\partial_{ij}\bar\varphi\partial_j\varphi
\en
thus
\be
T_6=\int -(\nabla\cdot\bar v)\Delta\bar\varphi\;\varphi + T_8 -T_3
\en
with $T_8=\int\partial_j\bar v_i\partial_{ij}\bar\varphi\varphi $.
Then 
\be
T_8 = &&-\int \partial_j\bar v_i\partial_j\bar\varphi\partial_i\varphi + \partial_{ij}\bar v_i\partial_j\bar\varphi\varphi\\
&&=-T_7+ \int \nabla\cdot\bar v(\Delta \bar\varphi\varphi+\nabla\bar\varphi \nabla\varphi)
\en
and finally we obtain
\be
S(t)=&&\int \nabla\bar v (\nabla\bar\varphi-\nabla\varphi)(\nabla\bar\varphi-\nabla\varphi) 
-\demi(\nabla\cdot \bar v)\left| \nabla\bar\varphi-\nabla\varphi\right|^2\\
&&+\demi\int(\nabla\cdot \bar v)\left| \nabla\varphi\right|^2 
\en
and the identity (\ref{1dtH3}) is proved.

$\hfill \Box$
\eject
\newpage


\section{Appendix}
\vspace{1cm}

\subsection*{Existence and uniqueness of solutions to second order ODE's with BV field}
In this section we prove existence and a.e. uniqueness for ordinary differential 
equations of   the form:
\beq\label{1ODE}
\Dt\left(\begin{array}{c} X\\ V\end{array}\right)=\left(\begin{array}{c} V \\ E(t,X)\end{array}\right)
\enq
for $X\in\Td,\;Y\in \Td$, and where the field $E$ belongs to $\Linf(]0,T[\times\Td)\cap L^1(]0,T[,BV(\Td))$. We work in the flat torus for simplicity, but our results are still
valid in an open subset of $\Rd$.
This result is an adaptation of the proof of \cite{DL} that uses the result 
of \cite{Bo} on renormalized solutions  of  transport equations.
\\
Remark: After this proof was written, the authors learned of a result by L. Ambrosio (\cite{Am}) that extends the results of \cite{DL} to  transport equations when the vector field is in $BV$ with bounded divergence.

\subsubsection*{Renormalized solutions for Vlasov equations with BV field}
Theorem 3.4  in \cite{Bo} adapted to the periodic case sates that if 
$f \in L^{\infty}(]0,T[\times \Td \times\Rd)$ satisfies:
\begin{eqnarray}\label{1kinebouchut}
\frac{\partial f}{\partial t}+\nabla_x\cdot\left(\xi f\right)
+\nabla_{\xi}\cdot\left( E(t,x)f\right)=0,
\end{eqnarray}
with  $E(t,x)\in L^1(]0,T[\times\Td)\cap L^{1}(]0,T[,BV(\Td))$, then for all $g$ Lipschitz continuous we have
\begin{eqnarray*}
\frac{\partial g(f)}{\partial t}+\nabla_x\cdot\left(\xi g(f)\right)
+\nabla_{\xi}\cdot\left( E(t,x)g(f)\right)=0.
\end{eqnarray*}
The property of renormalization implies that 
\begin{itemize}
\item solutions to (\ref{1kinebouchut})
with initial data in $L^{\infty}_{loc}(\Td \times\Rd)$ belong to 
\\$C(]0,T[, L^p_{loc}(\Td \times\Rd))$ for any $1\leq p <\infty$, 
\item solutions to (\ref{1kinebouchut}) with prescribed initial data in 
$L^{\infty}(\Td \times\Rd)$ are a.e. unique,
\item if $E_n$ converges to $E$ in $L^1(]0,T[\times \Td)$ then the solutions of (\ref{1kinebouchut}) with $E_n$ instead of $E$ converge to the solution of (\ref{1kinebouchut}).
\end{itemize} 
We notice that equation (\ref{1principaleper}) satisfies the assumptions of the Theorem, 
and thus will have the renormalization property. This renormalization 
property was used in \cite{DL} to obtain a.e. uniqueness for solutions 
of the corresponding ODE's. Indeed, the ODE
\be
&&\dt X(t,s,x)=b(t,X)\\
&&X(s,s,x)=x
\en
is associated to the transport equation:
\be
\dt u +b(t,x).\nabla u=0
\en
whose solutions satisfy for all $(t,s)\in]0,T[$
\be
u(t,X(t,s,x))=u(s,x).
\en
We extend this consequence to the case of second order equations, 
with $BV$ acceleration field.
To the kinetic equation
\beq\label{1cineapp}
\dt f +\xi\cdot\nabla_x f + E(t,x)\cdot\nabla_{\xi} f=0
\enq
we associate the second order ODE (\ref{1ODE}) which can rewritten as
 $\partial_{tt}{X}=E(t,X)$. The result is then the following:

\begin{theo}\label{1theo9}
Let $E(t,x) \in \Linf(]0,T[\times \Td)\cap L^1(]0,T[,BV(\Td))$,\\
 then the ODE
\begin{eqnarray}\label{1ODE2}
&&\partial_{tt}{X}(t,s,x,\xi)=E(t,X)\\
&&(X(s,s,x,\xi),\partial_t{X}(s,s,x,\xi))=(x,\xi)
\end{eqnarray}
admits an a.e. unique solution.
\end{theo}
Remark: Here almost everywhere must be understood for the Lebesgue measure of $\R^6$.
\\
{\bf Proof of Theorem \ref{1theo9}:}
We know that through equation (\ref{1ODE}) equation (\ref{1ODE2})
 can be considered as a first order differential equation. 
Let us first consider the 
case where $E$ is smooth. Note  $Y\in \Td\times\Rd $ (resp. $y$) for $(X,V)$ (resp. for $(x, \xi)$) 
and $B\in \Rd\times\Rd$ for $(\xi,E)$. Then for all $s\in ]0,T[$, $Y$ solves: 
\begin{eqnarray}
&&\dt Y(t,s,y)=B(t,Y(t,s,y))\label{1edo}\\
&&Y(s,s,y)=y\label{1edoin}
\end{eqnarray}
Then for all $t,t_1,t_2,t_3\in ]0,T[$ we have the following:
\begin{eqnarray*}
&&Y(t_3,t_2,Y(t_2,t_1,y))=Y(t_3,t_1,y)\\
&&Y(t,t,y)=y\\
&&Y(t_1,t_2,Y(t_2,t_1,y))=y.
\end{eqnarray*}
Differentiating the last equation with respect to $t_2$ yields:
\begin{eqnarray}
&&\partial_{s}Y(t,s,y)+\nabla_y Y(t,s,y)\cdot B(s,y)=0\label{1transport}\\
&&Y(t,t,y)=y\label{1transportin}.
\end{eqnarray}
$Y_t(s,y)=Y(t,s,y)$ thus solves a transport equation which is nothing but equation (\ref{1cineapp}). Using Theorem \ref{1Bouchut} we know that for all $g:\R^{2d}\rightarrow\R$ Lipschitz continuous,  $g(t,s,y)=g_0(Y(t,s,y))$ is the unique solution of 
\begin{eqnarray}
&&\partial_{s}g(t,s,y)+\nabla_y g(t,s,y)\cdot B(s,y)=0\label{1transg}\\
&&g(t,t,y)=g_0(y)\label{1transgin}.
\end{eqnarray}
Now we show existence and uniqueness for solutions of (\ref{1edo},\ref{1edoin}).
Let  $t$ and $s$ be fixed. Let us consider a regularization $E_n$ of the the field $E$ and set 
$B_n=(\xi,E_n)$. We have 
\begin{itemize}
\item$t\rightarrow Y_{1,n}(t,s,y)$ that satisfies (\ref{1edo},\ref{1edoin}) 
\item$s\rightarrow Y_{2,n}(t,s,y)$ that satisfies (\ref{1transport},\ref{1transportin}). 
\end{itemize}
From the stability Theorem 2.4 in \cite{DL} we know that  the whole sequence \\$t\rightarrow Y_{2,n}(t,s,.)$ converges in $C(]0,T[,L^p_{loc}(\Rd\times\Td))$ to $t\rightarrow Y_{2}(t,s,.)$, the unique renormalized solution of (\ref{1transport},\ref{1transportin}). Thus for fixed $t$ the whole sequence $Y_{2,n}(t,s,.)$ converges strongly in $L^p_{loc}(\Rd\times\Td)$. Now since for every $n$ we have  $Y_{1,n}(t,s,y)=Y_{2,n}(t,s,y)$ 
 the same property holds for $Y_{1,n}(s,t,.).$ Now we can pass to the limit in the term $B_n(t,Y_{1,n}(t,s,y))$. Indeed, by density of $C_c^{\infty}$ functions in $L^1,$ if we have $E_s\in C_c^{\infty}$ approximating $E$ then
\be
&&\|B(t,Y_n(t,s,y))-B(t,Y(t,s,y))\|_{L^1}\\
&\leq&\|B(t,Y_n(t,s,y))-B_s(t,Y_n(t,s,y))\|_{L^1}\\
&+&\|B_s(t,Y_n(t,s,y))-B_s(t,Y(t,s,y))\|_{L^1}\\
&+&\|B(t,Y(t,s,y))-B_s(t,Y(t,s,y))\|_{L^1}
\en
The second term goes to 0 because of the strong convergence of $Y_n$, the first and the third go to 0 because $Y$ and $Y_n$ are measure preserving mappings, and so for example $\|B(t,Y(t,s,y))-B_s(t,Y(t,s,y))\|_{L^1}=\|B(t,y)-B_s(t,y)\|_{L^1}$.
So finally we have 
\be
&&\|B_n(t,Y_n(t,s,y))-B(t,Y(t,s,y))\|_{L^1}\\
\leq&&\|B_n(t,Y_n(t,s,y))-B(t,Y_n(t,s,y))\|_{L^1}\\
&&+\|B(t,Y_n(t,s,y))-B(t,Y(t,s,y))\|_{L^1}
\en
that goes to 0 and we can pass to the limit in the equation (\ref{1edo},\ref{1edoin}) and the existence of a solution 
to (\ref{1edo},\ref{1edoin}) is proved. 

\medskip
\noindent
To obtain uniqueness, we argue as in \cite{DL}. Any function of the form $g_0(Y(t,s,y))$ is a solution of 
(\ref{1transg},\ref{1transgin}),  
thus by uniqueness of the solution of the transport equation we obtain uniqueness of the ODE.

$\hfill\Box$

\subsubsection*{A remark on ODE's of second order}
In this section, we want to solve the Cauchy problem for:
\begin{eqnarray*}
&&\partial_{tt}{X}(t,x)=E(t,X)\\
&&(X(0,x)\; ,\partial_t{X}(0,x))=(x,v(x))
\end{eqnarray*}
with $E$ as above. We are thus interested in monokinetic initial data.

\begin{theo}
for all $v^0(x)$ vector field on $\Td$, and for Lebesgue almost every $\delta v\in\mathbb R^d$, there exists an a.e.  unique solution to  
\begin{eqnarray*}
&&\partial_{tt}{X}(t,x)=E(t,X(t,x))\\
&&(X(0,x)\;,\partial_t{X}(0,x))=(x,v^0(x)+\delta v)
\end{eqnarray*}
\end{theo}

\medskip
\noindent
{\bf Proof:} Let $g(x,\xi)$ be the indicator function of the set of those $(x;\xi)$ such that the trajectory coming from x is not well defined. We just have to prove that for a.e. $\delta v\in \mathbb R^d$ we have $\int g(x,v^0(x)+\delta v)dx=0,$ which is true because
$$\int g(x,v^0(x)+\xi)dx d\xi=\int g(x,\xi)dx d\xi=0.$$

\paragraph{Stability}
Using the fact that for $E_n$ converging to $E$ in $L^1$ with\\
 $E\in L^1(]0,T[,BV(\Td))$, we have $X_n(t,x,v)\rightarrow X(t,x,v)$
 in $C([0,T],L^p)$, we have then, for all t, for almost every $\delta v,$ $X_n(t,x,v^0(x)+\delta v) \rightarrow X(t,x,v^0(x)+\delta v)$ in $L^p$. Thus we have 

\begin{theo}
If $E_n$ converges to $E$ in $L^1$ let  $X_n$ be solution of
\begin{eqnarray*}
&&\partial_{tt}{X_n}(t,x)=E_n(t,X_n(t,x))\\
&&(X_n(0,x),\partial_t{X_n}(0,x))=(x,v^0(x)+\delta v)
\end{eqnarray*}
 then for all t, for almost every  $\delta v,$ $X_n$ converges in $L^p(\mathbb R^3)-s$ to a solution (unique for almost every $\delta v$) of 
\begin{eqnarray*}
&&\partial_{tt}{X}(t,x)=E(t,X)\\
&&(X(0,x),\partial_t{X}(0,x))=(x,v^0(x)+\delta v).
\end{eqnarray*}
\end{theo}

\subsection*{Control of macroscopic density in kinetic\\ equations}
We prove here Lemma \ref{1densite}:
\begin{lemme}\label{1densite2}
Let $f\in \Linf([0,T]\times\Td\times\Rd)$ satisfy
\begin{eqnarray}
&&\frac{\partial f}{\partial t}+\nabla_x\cdot\left(\xi f\right)
+\nabla_{\xi}\cdot\left(E(t,x) f\right)=0\label{1cinetik1ap}\\
&&f(0,.,.)=f^0\label{1cinetik2ap}
\end{eqnarray}
with $E\in L^1([0,T]; BV(\Td))$ and
\beq
\|E\|_{\Linf([0,T]\times\Td)} \leq F.\label{1Ebv}
\enq
Let the initial condition $f_0$ be such that:
$$a(x,\xi)\leq f(0,x,\xi)\leq b(x,\xi),$$
with $\rho_a(x)=\int a(x,\xi)d\xi\geq m>0$ and $\rho_b(x)=\int b(x,\xi)d\xi\leq M<\infty$ and $a,b$ satisfying 
\beq
|\nabla_{x,\xi}a,b|\leq \frac{c}{1+|\xi|^{d+2}}\label{1lipab}.
\enq
Then there exists a constant $R>0$ such that
\be
(\rho_a(x)-Rt)\leq\rho(t,x)\leq(\rho_b(x)+Rt).
\en
\end{lemme}
{\bf Proof: } 
First suppose that the force field and the initial data are smooth. For equation (\ref{1cinetik1ap},\ref{1cinetik2ap})  we can exhibit characteristics  $(x,\xi)(t;t_0,x_0,\xi_0),$ giving the evolution of the particles in the phase space. We have $f(t,x,\xi)=f(t_0,x_0,\xi_0).$ Since the initial data is compactly supported and the force field is bounded in the $\Linf$ norm, we have 
$$|\xi-\xi_0|\leq F|t-t_0|,$$
$$|x-x_0|\leq (|\xi_0|+\frac{F}{2}|t-t_0|)|t-t_0|.$$

If for $t=0$ we have $a(x,\xi)\leq f(0,x,\xi)\leq b(x,\xi)$ then 
\begin{eqnarray*}
&&\underline{A}(t,x,\xi)\leq f(t,x,\xi)\leq\overline{B}(t,x,\xi)\\
&&\underline{A}(t,x,\xi)=\inf_{|\sigma_1|,|\sigma_2|\leq 1}a(x+|t-t_0|(\xi +\frac{F}{2}|t-t_0|)\sigma_1, \xi+F|t-t_0|\sigma_2)\\
&&\overline{B}(t,x,\xi)=\sup_{|\sigma_1|,|\sigma_2|\leq 1}b(x+|t-t_0|(\xi +\frac{F}{2}|t-t_0|)\sigma_1, \xi+F|t-t_0|\sigma_2).
\end{eqnarray*}
Using (\ref{1lipab}) and integrating in $\xi$ we find thus a constant $R=R(F,C,d)$ such that for $t-t_0\leq 1$ we have: 
\be
\rho_a(x)-R |t-t_0|\leq \rho(t,x)\leq \rho_b(x)+R|t-t_0|.
\en

\medskip
\noindent
Next we need to show that the solution of the regularized equation converges to the solution we are studying: this result comes from the uniqueness of the solution to (\ref{1cinetik1ap},\ref{1cinetik2ap}) which is a consequence of the renormalization property. Indeed since  $E$ is bounded in $BV$ the system (\ref{1cinetik1ap},\ref{1cinetik2ap}) admits a unique renormalized solution and the sequence of approximate solutions converge in $C([0,T],L^p_{x,\xi})$ for $1\leq p < \infty$ thus the bounds obtained above are preserved.

$\hfill$ $\Box$

\subsection*{Regularity of the polar factorization on the flat torus}
Here we deduce from \cite{Mc2}, \cite{Co} and \cite{Ca1}, \cite{Ca2}, \cite{Ca3} 
the Theorem \ref{1regper}. 
\begin{theo}\label{1regperap}
If $\rho \in C^{\alpha}(\Td)$ with $0<m\leq\rho\leq M$ is a probability measure on $\Td$ then $\Psi=\Psi[\rho]$ (see Definition \ref{1MaTdrho}) is a classical solution of 
\begin{eqnarray} 
\det D^2\Psi=\rho\label{1mongeperap}
\end{eqnarray}
and satisfies:
\begin{eqnarray} 
&&\|\nabla\Psi(x)-x\|_{\Linf}\leq C(d)=\sqrt d /2\\
&&\|D^2\Psi\|_{C^{\alpha}}\leq K(m,M,\|\rho\|_{C^{\alpha}})
\end{eqnarray}
\end{theo}
{\bf Proof of Theorem \ref{1regperap}: }Consider $\rho$ a $\Zd$ periodic probability measure, satisfying
\beq
0<m\leq\rho\leq M,\label{1mM}
\enq
and $\Phi[\rho]$ as in Definition \ref{1MaTdrho}.
First it is shown in \cite{Co} that 
\beq
|\nabla\Phi[\rho](x)-x|\leq C(d).\label{1LinfCo}
\enq
It follows that the strict convexity argument of \cite{Ca1} applies: indeed
if $\Phi=\Phi[\rho]$ is not strictly convex its graph contains a line and this 
contradicts (\ref{1LinfCo}). Moreover since $\Phi-|x|^2/2$ is globally Lipschitz and periodic there exists $N(d)$ such that $\|\Phi-|x|^2/2\|_{\Linf}\leq N(d)$. 
It follows then that there exists $0 <r(d)\leq R(d)$ and $M(d)$ such that 
\beq
B(r(d))\subset \{\Phi-\Phi(0) \leq M(d)\} \subset B(R(d))
\enq
It remains to show that our solution is a solution in the Aleksandrov sense
of the Monge-Amp\`ere equation 
\be
m\leq \det D^2 \Phi\leq M.
\en
This is a direct consequence with minor changes (to adapt to the periodic case) of Lemma 2 of \cite{Ca3}.
Then, normalizing $\Phi$ to $\tilde\Phi=\Phi-\Phi(0)-M(d)$ 
it follows that $\tilde\Phi$ is a solution of
\be
&&\rho(\nabla\tilde\Phi)\det D^2\tilde\Phi=1\\
&&\tilde\Phi = 0\;\;\text{ on }\partial\Omega\\
&&B(r(d))\subset \Omega \subset B(R(d))
\en
Thus the interior regularity results of \cite{Ca2} apply uniformly to
all $\Phi[\rho]$ with $\rho$ satisfying (\ref{1mM}) and $\|\rho\|_{C^{\alpha}(\Td)}$
bounded and Theorem \ref{1regper} follows.

$\hfill$ $\Box$

\subsection*{Acknowledgments}
The authors acknowledge the support
of the European IHP network "HYKE" HPRN-CT-2002-00282 and
the LRC CEA-Cadarache/UNSA.

\bibliography{vma-gafa-biblio}
\end{document}